# Self–Rectangulating Varieties of Type **5**

Keith A. Kearnes    Ágnes Szendrei [*]

**Abstract**

We show that a locally finite variety which omits abelian types is self–rectangulating if and only if it has a compatible semilattice term operation. Such varieties must have type–set {**5**}. These varieties are residually small and, when they are finitely generated, they have definable principal congruences. We show that idempotent varieties with a compatible semilattice term operation have the congruence extension property.

## 1 Introduction

In 1986, while investigating residually small varieties generated by a finite algebra, R. McKenzie discovered a property associated to nonabelian prime congruence quotients in members of such varieties which has become known as "rectangulation modulo a gene"; here a "gene" is a configuration of elements and polynomial operations associated to the nonabelian prime congruence quotient in question. Using this property McKenzie showed that whenever **A** is an algebra of size $n < \omega$ and $\mathcal{V} = \mathcal{V}(\mathbf{A})$ is residually small, then every finite subdirectly irreducible algebra in $\mathcal{V}$ which has nonabelian monolith has a tolerance $\tau$ with no more than $n$ tolerance blocks such that $\tau$ rectangulates itself modulo the gene associated with the monolith of the subdirectly irreducible. Since then McKenzie has used this result to characterize which locally finite varieties omitting abelian congruence quotients are residually small.

The rectangulation property is analyzed in [4] in a form where it is not localized to genes. There this property is shown to be a natural term condition which has a useful associated commutator operation. Finite "rectangular algebras", that is, finite algebras which are abelian with respect to this commutator, seem to be the nicest kinds of algebras among those which have a bound on the essential arity of their term operations. It turns out that this "global" rectangulation property is incompatible with the existence of genes, at least to the extent that no gene exists in a rectangular algebra. Hence "rectangulation modulo a gene" is a less restrictive notion than the rectangulation

[*]Research supported by the Hungarian National Foundation for Scientific Research grant no. T 17005

property from [4]. For binary relations $R$ and $S$ on an algebra, the fact that $R$ and $S$ rectangulate one another *modulo* a gene $G$ means roughly that, from the viewpoint of the gene $G$, the relations $R$ and $S$ "communicate with each other" in a rectangular way only.

In this paper we consider locally finite varieties of algebras where each finite member has the property that the total binary relation rectangulates itself modulo all genes associated to nonabelian prime congruence quotients. We will find that if in such a variety every nontrivial congruence quotient is nonabelian then the variety must have type–set $\{\mathbf{5}\}$. Therefore we call these varieties *self–rectangulating varieties of type–set* $\{\mathbf{5}\}$. The main result of this paper is an analogue of Herrmann's Theorem from modular commutator theory: we show that a locally finite variety of type–set $\{\mathbf{5}\}$ is self–rectangulating if and only if it has a semilattice term which commutes with all other term operations. The remainder of the paper is concerned with special properties of self–rectangulating varieties of type–set $\{\mathbf{5}\}$ (e.g., definable principal congruences, residual smallness, finite axiomatizability) with some emphasis on idempotent varieties.

## 2  Self–Rectangulation

The set of idempotent unary polynomial operations of an algebra $\mathbf{A}$ will be denoted by $E(\mathbf{A})$. If $\mathbf{A}$ is an algebra, $e \in E(\mathbf{A})$ and $1 \in e(A)$ is a fixed element, then we define a relation $\sqsubseteq$ on $A$ as follows. We say that $a \sqsubseteq b$ if and only if
$$ef(a) = 1 \Longrightarrow ef(b) = 1$$
for all $f \in \mathrm{Pol}_1(\mathbf{A})$. It is straightforward to check that $\sqsubseteq$ is a reflexive, transitive, compatible binary relation on $\mathbf{A}$, therefore it will be called the **natural quasiorder of A with respect to** $(e, 1)$.

It follows that for every natural quasiorder $\sqsubseteq$ as above, $\delta = \sqsubseteq \cap \sqsubseteq^{\cup}$ is a congruence of $\mathbf{A}$. Since $(a,b) \in \delta$ if and only if
$$ef(a) = 1 \Longleftrightarrow ef(b) = 1$$
for all $f \in \mathrm{Pol}_1(\mathbf{A})$, $\delta$ must be the largest congruence on $\mathbf{A}$ whose restriction $\delta|_{e(A)}$ to the set $e(A)$ has $\{1\}$ as a singleton class.

Now we give McKenzie's definition of "gene" and follow it with his definition of what it means for two relations to "rectangulate one another modulo a gene".

**Definition 2.1** Let $\mathbf{A}$ be an algebra and suppose that $0, 1 \in A$, $e(x) \in E(\mathbf{A})$, $x \sqcap y$ is a binary polynomial operation of $\mathbf{A}$ and $\rho$ is a congruence of $\mathbf{A}$. We say that $(e, \sqcap, 0, 1, \rho)$ is a **gene** of $\mathbf{A}$ if the following conditions hold:

(g1) $e(0) = 0 \neq 1 = e(1);$

(g2) $\rho$ is the largest congruence of $\mathbf{A}$ such that $\rho|_{e(A)}$ has $\{1\}$ as a singleton class;

(g3) $e(A)$ is closed under $\sqcap$, and for all $x \in e(A)$ we have $x \sqcap x = x$ and $x \sqcap 1 = x = 1 \sqcap x$;

(g4) if $a \in e(A) - \{1\}$, then $a$ is $\rho$–related to both $a \sqcap 0$ and $0 \sqcap a$.

There is a natural quasiorder corresponding to every gene $G = (e, \sqcap, 0, 1, \rho)$ of $\mathbf{A}$, namely the natural quasiorder $\sqsubseteq$ with respect to $(e, 1)$. Notice that in this case we have $\rho = \sqsubseteq \cap \sqsubseteq^\cup$, which is easily seen by comparing (g2) and the property of the congruence $\delta = \sqsubseteq \cap \sqsubseteq^\cup$ established right before Definition 2.1. Furthermore, the presence of the polynomial $\sqcap$ implies that within the set $e(A)$ the element 1 is maximal with respect to $\sqsubseteq$. Indeed, if $u \in e(A)$ and $u \neq 1$ then for the unary polynomial $f(x) = x \sqcap 1$ we have $ef(1) = e(1 \sqcap 1) = 1$ and $ef(u) = e(u \sqcap 1) = u$, hence $1 \not\sqsubseteq u$, as required.

Let $G = (e, \sqcap, 0, 1, \rho)$ and $G = (e', \sqcap', 0', 1', \rho')$ be arbitrary genes of an algebra $\mathbf{A}$. We will call $G$ and $G'$ **equivalent** if they have the same associated natural quasiorder. Clearly, if $G$ and $G'$ are equivalent then $\rho = \rho'$. An easy sufficient condition for $G$ and $G'$ to be equivalent is that $e(A) = e'(A)$ and $1 = 1'$; or more generally, that there is a polynomial isomorphism between $e(A)$ and $e'(A)$ which maps 1 to $1'$.

**Definition 2.2** Let $G = (e, \sqcap, 0, 1, \rho)$ be a gene of an algebra $\mathbf{A}$. Let $R, S$ be binary relations over $\mathbf{A}$. We say that $R$ **and** $S$ **rectangulate one another modulo** $G$ provided that whenever $p(\mathbf{x}, \mathbf{y}) \in \operatorname{Pol}_{m+n}(\mathbf{A})$, $\mathbf{a} R \mathbf{b}$ and $\mathbf{c} S \mathbf{d}$ we have
$$ep(\mathbf{a}, \mathbf{c}) = 1 = ep(\mathbf{b}, \mathbf{d}) \implies ep(\mathbf{a}, \mathbf{d}) = 1 = ep(\mathbf{b}, \mathbf{c}).$$

Later on, implications of the type displayed in the previous definition will usually be applied to equal relations $R = S$. Therefore for clarity we may underline arguments to indicate how we are applying the rule, as in
$$ep(\underline{\mathbf{a}}, \mathbf{c}) = 1 = ep(\mathbf{b}, \underline{\mathbf{d}}) \implies 1 = ep(\mathbf{a}, \mathbf{d}).$$

From the point of view of rectangulation equivalent genes can be considered the same. By this we mean that if the genes $G = (e, \sqcap, 0, 1, \rho)$ and $G' = (e', \sqcap', 0', 1', \rho')$ of an algebra $\mathbf{A}$ are equivalent then for any binary relations $R, S$ on $\mathbf{A}$, $R$ and $S$ rectangulate one another modulo $G$ if and only if they rectangulate one another modulo $G'$. To show this claim assume that $G$ and $G'$ have the same natural quasiorder $\sqsubseteq$ and that $R$ and $S$ rectangulate one another modulo $G'$. Let $p, \mathbf{a}, \mathbf{b}, \mathbf{c}, \mathbf{d}$ be as in Definition 2.2 and assume that $ep(\mathbf{a}, \mathbf{c}) = 1 = ep(\mathbf{b}, \mathbf{d})$. For any $f \in \operatorname{Pol}_1 \mathbf{A}$ with $e'f(1) = 1'$ the rectangulation property modulo $G'$, when applied to $e'fep(\mathbf{a}, \mathbf{c}) = 1' = e'fep(\mathbf{b}, \mathbf{d})$ yields that $e'fep(\mathbf{a}, \mathbf{d}) = 1' = e'fep(\mathbf{b}, \mathbf{c})$. Thus $1 \sqsubseteq ep(\mathbf{a}, \mathbf{d})$ and $1 \sqsubseteq ep(\mathbf{b}, \mathbf{c})$ with respect to $G'$. But all three elements appearing here belong to $e(A)$ and on this set 1 is maximal with respect to $\sqsubseteq$. Therefore $ep(\mathbf{a}, \mathbf{d}) = 1 = ep(\mathbf{b}, \mathbf{c})$. A similar argument with the roles of $G$ and $G'$ switched completes the proof.

It is a consequence of tame congruence theory that to every nonabelian prime congruence quotient $\langle \alpha, \beta \rangle$ on a finite algebra $\mathbf{A}$ there is an associated gene, called an $\langle \alpha, \beta \rangle$–**gene**, which is uniquely determined, up to equivalence. This gene is defined as follows. Let $U$ be an $\langle \alpha, \beta \rangle$–minimal set, let $e \in E(\mathbf{A})$ be such that $e(A) = U$, let $x \sqcap y$ be a binary polynomial of $\mathbf{A}$ which is a pseudo–meet operation of $\mathbf{A}|_U$, let $1 \in U$ be the unit element of this pseudo–meet operation, let $\rho$ be the pseudo–complement of $\beta$ over $\alpha$, and let $0 \in U - \{1\}$ be any element which is $\rho$–related to an element in the $\langle \alpha, \beta \rangle$–body of $U$. Then $(e, \sqcap, 0, 1, \rho)$ is a gene of $\mathbf{A}$. Furthermore, any two genes $(e, \sqcap, 0, 1, \rho)$ and $(e', \sqcap', 0', 1', \rho')$ constructed in this way are equivalent; in fact, there is a polynomial isomorphism between $e(A)$ and $e'(A)$ which sends 1 to $1'$. These claims can be proved from the material in Chapters 2, 4 and 5 of [2], especially the material on pages 29, 56–58 and 81. We emphasize that we shall only refer to an $\langle \alpha, \beta \rangle$–gene when $\langle \alpha, \beta \rangle$ is nonabelian.

The facts established in the preceding two paragraphs ensure that the following definition is meaningful.

**Definition 2.3** Let $\mathbf{A}$ be a finite algebra and $\langle \alpha, \beta \rangle$ a nonabelian prime congruence quotient of $\mathbf{A}$. We will say that **A rectangulates itself with respect to** $\langle \alpha, \beta \rangle$ if the total binary relation $T = A \times A$ rectangulates itself modulo an $\langle \alpha, \beta \rangle$–gene $G$. The natural quasiorder corresponding to $G$ will be called the **natural quasiorder of A with respect to** $\langle \alpha, \beta \rangle$. We will use the phrase **A is self–rectangulating** or **A rectangulates itself** to mean that $\mathbf{A}$ rectangulates itself with respect to every nonabelian prime congruence quotient. Finally, a locally finite variety will be called self–rectangulating if every finite member is.

In this paper we will mainly use self–rectangulation with respect to $\langle 0_\mathbf{A}, \mu \rangle$ where $\mathbf{A}$ is a subdirectly irreducible algebra with nonabelian monolith $\mu$. If $G = (e, \sqcap, 0, 1, \rho)$ is a $\langle 0_\mathbf{A}, \mu \rangle$–gene then $\rho = 0_\mathbf{A}$, and hence the corresponding natural quasiorder $\sqsubseteq$ is a partial order. Therefore in this case $\sqsubseteq$ will be called the **natural order of A with respect to** $\langle 0_\mathbf{A}, \mu \rangle$.

Now we explain why the property of being self–rectangulating forces every nonabelian prime quotient to be of type $\mathbf{5}$. Let $\mathbf{A}$ be a finite algebra, let $\langle \alpha, \beta \rangle$ be a nonabelian prime quotient of $\mathbf{A}$, and assume that $\mathrm{typ}(\alpha, \beta) \in \{\mathbf{3}, \mathbf{4}\}$. Let $G = (e, \sqcap, 0, 1, \rho)$ be an $\langle \alpha, \beta \rangle$–gene defined on $U = e(A)$, which is an $\langle \alpha, \beta \rangle$–minimal set. Choose a polynomial operation $\sqcup$ which is an $\langle \alpha, \beta \rangle$ pseudo–join operation for $\mathbf{A}|_U$. Then the following displayed line proves that $\mathbf{A}$ does not rectangulate itself with respect to $\langle \alpha, \beta \rangle$:

$$e(\underline{0} \sqcup 1) = 1 = e(1 \sqcup \underline{0}) \qquad \text{but} \qquad e(\underline{0} \sqcup \underline{0}) = e(0) = 0.$$

We conclude that if all prime quotients of $\mathbf{A}$ are nonabelian and $\mathbf{A}$ is self–rectangulating, then $\mathrm{typ}\{\mathbf{A}\} = \{\mathbf{5}\}$.

Note that if **A** rectangulates itself with respect to a nonabelian prime quotient $\langle \alpha, \beta \rangle$ and $\delta \leq \alpha$, then $\mathbf{A}/\delta$ rectangulates itself with respect to $\langle \alpha/\delta, \beta/\delta \rangle$. This is a consequence of the fact that the $\delta|_U$–class of the element 1 is just $\{1\}$ when $\delta \leq \rho$. It follows that when a finite algebra **A** is self–rectangulating, then every homomorphic image of **A** is also self–rectangulating.

**Definition 2.4** If **B** is an algebra, then a **compatible semilattice operation** of **B** is an idempotent, commutative, associative binary operation $\wedge$ on $B$ which is a homomorphism $\wedge : \mathbf{B}^2 \to \mathbf{B}$. If $\wedge$ agrees with a term [polynomial] operation of **B**, we call the term [polynomial] operation a **compatible semilattice term [polynomial] operation** of **B**. A compatible semilattice term operation of a variety $\mathcal{V}$ is a term operation which interprets as a compatible semilattice operation in every member of $\mathcal{V}$.

Note that a compatible semilattice operation $\wedge$ of **B** is compatible with all polynomial operations of **B** as well, since $\wedge$ is idempotent.

**LEMMA 2.5** *If an algebra* **B** *has a compatible semilattice operation* $\wedge$*, then every semilattice polynomial operation of* **B** *coincides with* $\wedge$.

*Proof.* Let $p$ be any binary, idempotent, commutative polynomial operation of **B**. Then $p$ commutes with $\wedge$, so for any elements $a, b \in B$ we have

$$\begin{aligned} p(a,b) &= p(a,b) \wedge p(a,b) \\ &= p(a,b) \wedge p(b,a) \\ &= p(a \wedge b, b \wedge a) \\ &= p(a \wedge b, a \wedge b) \\ &= a \wedge b. \end{aligned}$$

This proves the claim. $\square$

The following corollary is an immediate consequence.

**COROLLARY 2.6** *An algebra has at most one compatible semilattice polynomial operation.*

The next lemma describes a construction due to McKenzie, which is an important step in finding a characterization for self–rectangulation. Let **A** be an algebra and $N$ any positive integer. In the algebra $\mathbf{A}^N$ we let $\widehat{a}$ denote the constant function with range $\{a\}$ for any $a \in A$. For any polynomial $p \in \mathrm{Pol}\,(\mathbf{A})$ we let $\widehat{p}$ denote the polynomial of $\mathbf{A}^N$ which is $p$ acting coordinatewise.

**LEMMA 2.7** *Let* **A** *be a finite subdirectly irreducible algebra with nonabelian monolith* $\mu$. *Consider a* $\langle 0_{\mathbf{A}}, \mu \rangle$*–gene* $G = (e, \sqcap, 0, 1, 0_{\mathbf{A}})$ *where* $U = e(A)$ *is a* $\langle 0_{\mathbf{A}}, \mu \rangle$*–minimal set. Furthermore, let* $N$ *be a positive integer, let* $\sqsubseteq$ *be the natural quasiorder of* $\mathbf{A}^N$ *with respect to* $(\widehat{e}, \widehat{1})$, *and let* $\delta = \sqsubseteq \cap \sqsubseteq^{\cup}$. *If* **A** *rectangulates itself with respect to* $\langle 0_{\mathbf{A}}, \mu \rangle$, *then the algebra* $\mathbf{A}^N/\delta$ *has the following properties:*

(1) **A** can be embedded in $\mathbf{A}^N/\delta$;

(2) $\mathbf{A}^N/\delta$ is subdirectly irreducible; and

(3) for sufficiently large $N$ the algebra $\mathbf{A}^N/\delta$ has a compatible semilattice operation.

*Proof.* It follows from the construction of $G$ that the body of $U$ is the set $\{0, 1\}$. We may assume that $e(x \sqcap y) = x \sqcap y$.

First we prove (2). Since $\delta$ is the largest congruence on $\mathbf{A}^N$ whose restriction to the set $\widehat{e}(A^N) = U^N$ has $\{\widehat{1}\}$ as a singleton class, we know that $(\widehat{0}, \widehat{1}) \notin \delta$. We claim that $\delta^* := \mathrm{Cg}(\delta \cup \{(\widehat{0}, \widehat{1})\})$ is the unique upper cover of $\delta$.

Choose any $\mathbf{u}, \mathbf{v} \in A^N$ such that $(\mathbf{u}, \mathbf{v}) \notin \delta$. Then there is an $f \in \mathrm{Pol}_1(\mathbf{A}^N)$ such that, say, $\widehat{e}f(\mathbf{u}) = \widehat{1}$ and $\widehat{e}f(\mathbf{v}) = \mathbf{r} \neq \widehat{1}$. We shall argue that $(\mathbf{r}, \mathbf{r} \sqcap \widehat{0}) \in \delta$. Then, since $(\widehat{1}, \mathbf{r}) \in \mathrm{Cg}(\mathbf{u}, \mathbf{v})$, and therefore $(\widehat{1} \sqcap \widehat{0}, \mathbf{r} \sqcap \widehat{0}) \in \mathrm{Cg}(\mathbf{u}, \mathbf{v})$, we will have that

$$\widehat{1} \; \mathrm{Cg}(\mathbf{u}, \mathbf{v}) \; \mathbf{r} \; \delta \; \mathbf{r} \sqcap \widehat{0} \; \mathrm{Cg}(\mathbf{u}, \mathbf{v}) \; \widehat{1} \sqcap \widehat{0} = \widehat{0}$$

which tells us that $(\widehat{1}, \widehat{0}) \in \delta \vee \mathrm{Cg}(\mathbf{u}, \mathbf{v})$ whenever $(\mathbf{u}, \mathbf{v}) \notin \delta$. This will prove that $\delta^*$ is not contained in $\delta$ but is contained in any congruence larger that $\delta$.

Since $\mathbf{r} \in U^N - \{\widehat{1}\}$, we get that there is some coordinate $i$ such that $r_i \in U - \{1\}$. Since $x \sqcap 0 = x$ for $x \in U - \{1\}$, we get that $r_i = r_i \sqcap 0 = (\mathbf{r} \sqcap \widehat{0})_i$. That is, $\mathbf{r}$ and $\mathbf{r} \sqcap \widehat{0}$ agree in the $i$–th coordinate. Assume that for some $f \in \mathrm{Pol}_1(\mathbf{A}^N)$ we have $\widehat{e}f(\mathbf{r}) = \widehat{1}$. We claim that $\widehat{e}f(\mathbf{r} \sqcap \widehat{0}) = \widehat{1}$ as well. To see this, first write $f$ as

$$f(x) = \widehat{g}(x, \mathbf{w}^1, \ldots, \mathbf{w}^k)$$

for some $g \in \mathrm{Pol}(\mathbf{A})$ and tuples $\mathbf{w}^l \in A^N$ ($1 \leq l \leq k$). Let $\mathbf{w}_i$ denote the tuple of $i$–th coordinates of the tuples $\mathbf{w}^1, \ldots, \mathbf{w}^k$. Thus in the $i$–th coordinate we have

$$eg(r_i \sqcap 0, \mathbf{w}_i) = eg(r_i, \mathbf{w}_i) = (\widehat{e}f(\mathbf{r}))_i = (\widehat{1})_i = 1.$$

If $j$ is any coordinate whatsoever, we have

$$eg(r_j, \mathbf{w}_j) = (\widehat{e}f(\mathbf{r}))_j = (\widehat{1})_j = 1.$$

Therefore, using the fact that $r_j \sqcap 1 = r_j$, we may apply the rectangulation condition to

$$eg(r_i \sqcap \underline{0}, \mathbf{w}_i) = 1 = eg(r_j, \mathbf{w}_j) = eg(\underline{r_j} \sqcap 1, \underline{\mathbf{w}_j})$$

to obtain that $eg(r_j \sqcap 0, \mathbf{w}_j) = 1$ holds for all $j$. Hence $\widehat{e}f(\mathbf{r} \sqcap \widehat{0}) = \widehat{1}$ as we claimed. Therefore,

$$\widehat{e}f(\mathbf{r}) = \widehat{1} \implies \widehat{e}f(\mathbf{r} \sqcap \widehat{0}) = \widehat{1}$$

and so $\mathbf{r} \sqsubseteq \mathbf{r} \sqcap \widehat{0}$ in $\mathbf{A}^N$. The reverse inclusion is proved in exactly the same fashion. This completes the proof that $(\mathbf{r}, \mathbf{r} \sqcap \widehat{0}) \in \delta$ and therefore that $\delta^*$ is the unique upper cover of $\delta$. Thus $\mathbf{B} = \mathbf{A}^N/\delta$ is a subdirectly irreducible algebra in $\mathcal{V}(\mathbf{A})$.

To verify (1) observe that the composition of the homomorphism $\Delta\colon \mathbf{A} \to \mathbf{A}^N$, $x \mapsto \widehat{x}$ followed by the homomorphism $\eta\colon \mathbf{A}^N \to \mathbf{B}$, $x \mapsto x/\delta$ is a homomorphism from $\mathbf{A}$ to $\mathbf{B}$. The monolith of $\mathbf{A}$ is $\mathrm{Cg}(0,1)$ and

$$(\eta \circ \Delta)(0) = \widehat{0}/\delta \neq \widehat{1}/\delta = (\eta \circ \Delta)(1).$$

Hence $\eta \circ \Delta$ is one–to–one and therefore an isomorphism of $\mathbf{A}$ onto a subalgebra of $\mathbf{B}$.

What remains to show is that for sufficiently large $N$, $\mathbf{B}$ has a compatible semilattice operation. Let $N \geq 4|A|$. For each $\mathbf{a} = (a_1, \ldots, a_N) \in A^N$ we define a subset $\mathrm{ran}(\mathbf{a}) \subseteq A$ as follows:

$$\mathrm{ran}((a_1, \ldots, a_N)) = \{a_1, \ldots, a_N\}.$$

We argue now that if $\mathrm{ran}(\mathbf{a}) \supseteq \mathrm{ran}(\mathbf{b})$, then $\mathbf{a} \sqsubseteq \mathbf{b}$ in $\mathbf{A}^N$. Assume that $\mathrm{ran}(\mathbf{a}) \supseteq \mathrm{ran}(\mathbf{b})$ and that for some $f \in \mathrm{Pol}_1(\mathbf{A}^N)$ we have $\widehat{e}f(\mathbf{a}) = \widehat{1}$. As above, we may write $f$ as

$$f(x) = \widehat{g}(x, \mathbf{w}^1, \ldots, \mathbf{w}^k)$$

for some $g \in \mathrm{Pol}(\mathbf{A})$ and tuples $\mathbf{w}^l \in A^N$ ($1 \leq l \leq k$). Our goal is to prove that

$$\widehat{1} = \widehat{e}f(\mathbf{b}) = \widehat{e}\widetilde{g}(\mathbf{b}, \mathbf{w}^1, \ldots, \mathbf{w}^k).$$

It suffices to show that in the $i$–th coordinate we have $eg(b_i, \mathbf{w}_i) = 1$ for an arbitrarily chosen $i$. For any coordinate $i$, we have $b_i \in \mathrm{ran}(\mathbf{b}) \subseteq \mathrm{ran}(\mathbf{a})$, so there is a $j$ with $a_j = b_i$. Since $\widehat{e}f(\mathbf{a}) = \widehat{1}$, we have in the $i$–th and $j$–th coordinates that

$$eg(a_i, \underline{\mathbf{w}_i}) = 1 = eg(a_j, \mathbf{w}_j) = eg(\underline{b_i}, \mathbf{w}_j).$$

Using the rectangulation rule we get that $eg(b_i, \mathbf{w}_i) = 1$ as we hoped. This completes the proof that $\mathbf{a} \sqsubseteq \mathbf{b}$. Note that this implies that if $\mathrm{ran}(\mathbf{a}) = \mathrm{ran}(\mathbf{b})$, then $(\mathbf{a}, \mathbf{b}) \in \delta$.

We describe an operation $\wedge$ on $B$ according to the following rule: For $\mathbf{a}/\delta, \mathbf{b}/\delta \in B$ we define

$$(\mathbf{a}/\delta) \wedge (\mathbf{b}/\delta) = (\mathbf{c}/\delta)$$

where $\mathbf{c} \in A^N$ is any tuple for which $\mathrm{ran}(\mathbf{c}) = \mathrm{ran}(\mathbf{a}) \cup \mathrm{ran}(\mathbf{b})$. There are plenty of tuples $\mathbf{c}$ with this property since the number $N \geq 4|A|$ of coordinates exceeds $|A|$. Our first task is to show that $\wedge$ is well defined. Assume that $(\mathbf{a}, \mathbf{a}'), (\mathbf{b}, \mathbf{b}') \in \delta$ and that $\mathrm{ran}(\mathbf{c}) = \mathrm{ran}(\mathbf{a}) \cup \mathrm{ran}(\mathbf{b})$ and $\mathrm{ran}(\mathbf{c}') = \mathrm{ran}(\mathbf{a}') \cup \mathrm{ran}(\mathbf{b}')$. We shall argue that $\mathbf{c} \sqsubseteq \mathbf{c}'$. Coupling this with an identical argument with the roles of $\mathbf{c}$ and $\mathbf{c}'$ switched proves that $(\mathbf{c}, \mathbf{c}') \in \delta$, so $\wedge$ is well defined. Since $\mathrm{ran}(\mathbf{c}) \supseteq \mathrm{ran}(\mathbf{a}), \mathrm{ran}(\mathbf{b})$ the result of the last paragraph shows that $\mathbf{c} \sqsubseteq \mathbf{a}, \mathbf{b}$. Pick any $f \in \mathrm{Pol}_1(\mathbf{A}^N)$ such that $\widehat{e}f(\mathbf{c}) = \widehat{1}$. Since $\mathbf{c} \sqsubseteq \mathbf{a}, \mathbf{b}$ we get that $\widehat{e}f(\mathbf{a}) = \widehat{e}f(\mathbf{b}) = \widehat{1}$, too. Since $(\mathbf{a}, \mathbf{a}'), (\mathbf{b}, \mathbf{b}') \in \delta$ we get that $\widehat{e}f(\mathbf{a}') = \widehat{e}f(\mathbf{b}') = \widehat{1}$. As we have done twice before, write $f(x)$ as $\widehat{g}(x, \mathbf{w}^1, \ldots, \mathbf{w}^k)$. To prove that $\widehat{e}f(\mathbf{c}') = \widehat{1}$, which is our current goal, it suffices to prove that for an arbitrarily chosen coordinate $i$ we have $eg(c_i', \mathbf{w}_i) = 1$. Since

$$c_i' \in \mathrm{ran}(\mathbf{c}') = \mathrm{ran}(\mathbf{a}') \cup \mathrm{ran}(\mathbf{b}')$$

there is some $a_j'$ or $b_j'$ equal to $c_i'$. Assuming the former, we have

$$eg(\underline{c_i'}, \mathbf{w}_j) = eg(a_j', \mathbf{w}_j) = 1 = eg(a_i', \underline{\mathbf{w}_i})$$

which implies that $eg(c_i', \mathbf{w}_i) = 1$. This finishes the proof that $\mathbf{c} \sqsubseteq \mathbf{c}'$ and also the proof that $\wedge$ is well defined.

That $\wedge$ is idempotent, commutative and associative follows from the fact that the operation of union has these properties.

What remains to show is that $\wedge$ is a homomorphism from $\mathbf{B}^2$ to $\mathbf{B}$. Let $t(x_1, \ldots, x_m)$ be a basic operation of $\mathbf{B}$. We must show that

$$t(\mathbf{a}_1/\delta, \ldots, \mathbf{a}_m/\delta) \wedge t(\mathbf{b}_1/\delta, \ldots, \mathbf{b}_m/\delta) = t(\mathbf{a}_1/\delta \wedge \mathbf{b}_1/\delta, \ldots, \mathbf{a}_m/\delta \wedge \mathbf{b}_m/\delta)$$

for any $\mathbf{a}_1, \ldots, \mathbf{a}_m, \mathbf{b}_1, \ldots, \mathbf{b}_m \in A^N$. Let $R = t(\mathrm{ran}(\mathbf{a}_1), \ldots, \mathrm{ran}(\mathbf{a}_m)) \subseteq A$ and let $S = t(\mathrm{ran}(\mathbf{b}_1), \ldots, \mathrm{ran}(\mathbf{b}_m)) \subseteq A$. Let $\Omega$ be the set of tuples in $A^N$ whose first $2|A|$ components are from $R$ and which exhaust $R$, and whose last $2|A|$ components are from $S$ and exhaust $S$. Since all $\mathbf{p} \in \Omega$ have $\mathrm{ran}(\mathbf{p}) = R \cup S$, all members of $\Omega$ are $\delta$–related. We claim first that $t(\mathbf{a}_1, \ldots, \mathbf{a}_m)$ is $\delta$–related to a tuple $\mathbf{a}$ with $\mathrm{ran}(\mathbf{a}) = R$ and that $t(\mathbf{b}_1, \ldots, \mathbf{b}_m)$ is $\delta$–related to a tuple $\mathbf{b}$ with $\mathrm{ran}(\mathbf{b}) = S$. Once we prove this claim it will establish that the left hand side of the previous displayed line equals $\mathbf{p}/\delta$ for some (any) $\mathbf{p} \in \Omega$. Afterwards we will prove that the right hand side has the same value. To show that $t(\mathbf{a}_1, \ldots, \mathbf{a}_m)$ is $\delta$–related to a tuple $\mathbf{a}$ with $\mathrm{ran}(\mathbf{a}) = R$, we shall argue that it is possible to alter each $\mathbf{a}_i$ to a $\delta$–related tuple $\mathbf{a}_i'$ such that $\mathrm{ran}(t(\mathbf{a}_1', \ldots, \mathbf{a}_m')) = R$. Enumerate $R$ as $(r_1, \ldots, r_l)$. For $j$ ranging from 1 to $l$, choose $a_{ij}' \in \mathrm{ran}(\mathbf{a}_i)$ so that $t(a_{1j}', \ldots, a_{lj}') = r_j$. For $j$ ranging from $l+1$ to $N$, choose any $a_{ij}' \in \mathrm{ran}(\mathbf{a}_i)$ so long as the choices for a fixed $i$ include a full list of all elements of $\mathrm{ran}(\mathbf{a}_i)$. There are enough coordinates to do this since $|A| \leq N - l$. We now have tuples $\mathbf{a}_i' = (a_{i1}', \ldots, a_{iN}')$ with $\mathrm{ran}(\mathbf{a}_i) = \mathrm{ran}(\mathbf{a}_i')$, so $(\mathbf{a}_i, \mathbf{a}_i') \in \delta$, such that

$$t(\mathbf{a}_1', \ldots, \mathbf{a}_m') = (r_1, \ldots, r_l, \ldots) =: \mathbf{a}.$$

The left hand side of this display is $\delta$–related to $t(\mathbf{a}_1, \ldots, \mathbf{a}_m)$ while $\mathrm{ran}(\mathbf{a}) = R$. This finishes the proof of our claim that

$$t(\mathbf{a}_1/\delta, \ldots, \mathbf{a}_m/\delta) \wedge t(\mathbf{b}_1/\delta, \ldots, \mathbf{b}_m/\delta) = \mathbf{p}/\delta$$

for some $\mathbf{p} \in \Omega$. Now we must show that

$$t(\mathbf{a}_1/\delta \wedge \mathbf{b}_1/\delta, \ldots, \mathbf{a}_m/\delta \wedge \mathbf{b}_m/\delta) = \mathbf{p}/\delta$$

for some $\mathbf{p} \in \Omega$. It is possible to construct $\mathbf{c}_1, \ldots, \mathbf{c}_m \in A^N$ with the following properties:

(a) Each $\mathbf{c}_i$ has its first $2|A|$ coordinates chosen so that they come from $\mathrm{ran}(\mathbf{a}_i)$ and exhaust $\mathrm{ran}(\mathbf{a}_i)$.

(b) Each $\mathbf{c}_i$ has its last $2|A|$ coordinates chosen so that they come from $\mathrm{ran}(\mathbf{b}_i)$ and exhaust $\mathrm{ran}(\mathbf{b}_i)$.

(c) The coordinates of the tuples $\mathbf{c}_i$ are arranged so that the first $2|A|$ coordinate values of the tuple $t(\mathbf{c}_1,\ldots,\mathbf{c}_m)$ exhaust $R$ and the last $2|A|$ values exhaust $S$.

The procedure for constructing the $\mathbf{c}_i$ is similar to the procedure used earlier in this paragraph. Since each $\mathrm{ran}(\mathbf{c}_i) = \mathrm{ran}(\mathbf{a}_i) \cup \mathrm{ran}(\mathbf{b}_i)$, we have that $\mathbf{a}_i/\delta \wedge \mathbf{b}_i/\delta = \mathbf{c}_i/\delta$. Thus

$$\begin{aligned} t(\mathbf{a}_1/\delta \wedge \mathbf{b}_1/\delta, \ldots, \mathbf{a}_m/\delta \wedge \mathbf{b}_m/\delta) &= t(\mathbf{c}_1/\delta, \ldots, \mathbf{c}_m/\delta) \\ &= t(\mathbf{c}_1, \ldots, \mathbf{c}_m)/\delta \\ &= \mathbf{p}/\delta \end{aligned}$$

where $\mathbf{p} \in \Omega$. The last equality holds since the first $2|A|$ coordinates of $t(\mathbf{c}_1,\ldots,\mathbf{c}_m)$ are contained in $R$ and exhaust $R$ while the last $2|A|$ coordinates are contained in $S$ and exhaust $S$. This establishes that $\wedge$ is compatible with $t$, which was arbitrarily chosen, so $\wedge$ is a homomorphism. $\square$

**THEOREM 2.8** *Let $\mathbf{A}$ be a finite algebra which omits types **1** and **2**. The following conditions are equivalent.*

(i) $\mathbf{A}$ *is self–rectangulating.*

(ii) $\mathbf{A}$ *is isomorphic to a subalgebra of a finite algebra $\mathbf{B} \in \mathcal{V}(\mathbf{A})$ such that $\mathbf{B}$ has a compatible semilattice operation.*

*If $\mathbf{A}$ is subdirectly irreducible and $\mu$ denotes its monolith, then condition* (iii) *below is also equivalent to* (i) *and* (ii):

(iii) $\mathbf{A}$ *rectangulates itself with respect to $\langle 0_\mathbf{A}, \mu \rangle$.*

*Proof.* Any finite $\mathbf{A}$ which omits types **1** and **2** and is self–rectangulating is a subdirect product of subdirectly irreducible algebras having all of these properties. Furthermore, the class of finite algebras $\mathbf{A}$ which have finite extensions in $\mathcal{V}(\mathbf{A})$ with compatible semilattice operations is closed under the formation of subalgebras and finite products. Therefore it suffices to prove (i) $\Rightarrow$ (ii) for the case when $\mathbf{A}$ is subdirectly irreducible. However, then (i) obviously implies (iii) and Lemma 2.7 proves that (iii) implies (ii).

That (ii) implies (i) can be verified as follows. Assume that $\mathbf{A}$ is a subalgebra of $\mathbf{B}$ where $\mathbf{B}$ has a compatible semilattice operation $\wedge$. Select a gene $(e, \sqcap, 0, 1, \rho)$ in $\mathbf{A}$, and let $U = e(A)$. We claim that the element 1 is maximal in $U$ under the restriction of the semilattice ordering $\leq$ of $\mathbf{B}$ to $U$. This follows by observing that if $u \in U$ is such that $1 \leq u$, then

$$1 = 1 \sqcap 1 = (u \wedge 1) \sqcap (1 \wedge u) = (u \sqcap 1) \wedge (1 \sqcap u) = u \wedge u = u.$$

To verify self–rectangulation choose a polynomial $f \in \operatorname{Pol}_{m+n}(\mathbf{A})$ and tuples $\mathbf{a}, \mathbf{b} \in A^m$, $\mathbf{c}, \mathbf{d} \in A^n$ such that

$$ef(\underline{\mathbf{a}}, \mathbf{c}) = 1 = ef(\mathbf{b}, \underline{\mathbf{d}}).$$

We have that $ef(\mathbf{a}, \mathbf{d}) \in U$, since $e(A) = U$, and

$$\begin{aligned}
1 &= 1 \wedge 1 \\
&= ef(\mathbf{a}, \mathbf{c}) \wedge ef(\mathbf{b}, \mathbf{d}) \\
&= ef(\mathbf{a} \wedge \mathbf{b}, \mathbf{c} \wedge \mathbf{d}) \\
&= ef(\mathbf{a} \wedge \mathbf{b}, \mathbf{d} \wedge \mathbf{c}) \\
&= ef(\mathbf{a}, \mathbf{d}) \wedge ef(\mathbf{b}, \mathbf{c})
\end{aligned}$$

which shows that $1 \leq ef(\mathbf{a}, \mathbf{d})$ in the restriction of the semilattice order to $U$. Since 1 is maximal, we deduce that $ef(\mathbf{a}, \mathbf{d}) = 1$, which establishes the rectangulation condition. This completes the proof. $\square$

**COROLLARY 2.9** *Let $\mathbf{A}$ be a finite algebra such that $\mathbf{A}$ omits types $\mathbf{1}$ and $\mathbf{2}$ and $\mathcal{V}(\mathbf{A})$ omits type $\mathbf{1}$. If $\mathbf{A}$ is self–rectangulating, then $\mathcal{V}(\mathbf{A})$ is self–rectangulating and of type–set $\{\mathbf{5}\}$.*

*Proof.* The previous theorem shows that if $\mathbf{A}$ is self–rectangulating, then $\mathbf{A}$ has an extension $\mathbf{B}$ which has a compatible semilattice operation. Let $\mathbf{F}$ be a finitely generated free algebra in $\mathcal{V}(\mathbf{A})$. We have that $\mathbf{F} \in \mathsf{SP}(\mathbf{A}) \subseteq \mathsf{SP}(\mathbf{B})$ and so $\mathbf{F}$ is isomorphic to a subalgebra of an algebra $\mathbf{G}$ (some power of $\mathbf{B}$) which has a compatible semilattice operation $\wedge$.

We show that $\mathbf{F}$ omits type $\mathbf{2}$. Suppose not, and let $U$ be a minimal set for a prime congruence quotient of type $\mathbf{2}$. Let $d'$ be a ternary polynomial of $\mathbf{F}$ which is a pseudo–Mal'cev operation for $\mathbf{F}|_U$. Clearly, $\mathbf{G}$ has a polynomial operation $d$ extending $d'$. Calculating in $\mathbf{G}$ we see that for any elements $u, v$ in the body of $U$ we have

$$\begin{aligned}
u &= d(u, v, v) \wedge d(v, v, u) \\
&= d(u \wedge v, v \wedge v, v \wedge u) \\
&= d(u \wedge v, v \wedge v, u \wedge v) \\
&= d(u, v, u) \wedge v,
\end{aligned}$$

whence $u \leq v$ in the semilattice order of $\mathbf{G}$. Switching the role of $u$ and $v$ we conclude that $u = v$. This implies that the body of $U$ is a singleton, which is impossible.

Combining the assumption on $\mathcal{V}(\mathbf{A})$ with the fact proved above we get that $\operatorname{typ}\{\mathbf{F}\} \cap \{\mathbf{1}, \mathbf{2}\} = \emptyset$. Therefore we can apply the previous theorem once again to deduce that **F** is self–rectangulating. Every finite algebra is a homomorphic image of a finitely generated, self–rectangulating, free algebra in $\mathcal{V}(\mathbf{A})$. This proves that every finite algebra in $\mathcal{V}(\mathbf{A})$ is self–rectangulating. Since $\mathcal{V}(\mathbf{A})$ omits types **1** and **2** and is self–rectangulating, it follows that $\mathrm{typ}\{\mathcal{V}(\mathbf{A})\} = \{\mathbf{5}\}$. □

In Corollary 2.9 the assumption that type **1** does not occur in $\mathcal{V}(\mathbf{A})$ cannot be omitted, as the next example shows.

**Example 2.10** Let $\mathbf{A} = (A; \cdot, e, f, 0)$ be the 3–element algebra where $A = \{0, 1, 2\}$, 0 denotes the unary constant operation with value 0 and the remaining operations are defined as follows:

$$x \cdot y = \begin{cases} 0 & \text{if } x = 0 \\ y & \text{otherwise} \end{cases},$$

and

$$e(x) = \begin{cases} 0 & \text{if } x = 0 \\ 2 & \text{otherwise} \end{cases}, \qquad f(x) = \begin{cases} 0 & \text{if } x \neq 2 \\ 1 & \text{if } x = 2 \end{cases}.$$

It is straightforward to check that **A** is a simple algebra of type **5**. The meet operation of the chain $0 \leq 1 \leq 2$ is easily seen to be a compatible semilattice operation of **A**, hence by Theorem 2.8 **A** is self–rectangulating. We show that the variety $\mathcal{V}(\mathbf{A})$ has exactly two subdirectly irreducible algebras: **A** and the 4–element algebra $\mathbf{S} = \mathbf{C}/\gamma$ where **C** is the subalgebra of $\mathbf{A}^2$ with universe $A^2 - \{(2,1)\}$ and $\gamma$ is the congruence collapsing all pairs having at least one 0 coordinate.

It suffices to verify that **A** and **S** are the only finite subdirectly irreducible algebras in $\mathcal{V}(\mathbf{A})$. The following identities hold in **A**: $x \cdot (y \cdot z) = (x \cdot y) \cdot (x \cdot z)$, $x \cdot e(y) = e(x \cdot y)$, $x \cdot f(y) = f(x \cdot y)$ and $x \cdot 0 = 0$. These equations imply that for any $\mathbf{B} \in \mathcal{V}(\mathbf{A})$ and any $b \in B$ the function $x \mapsto b \cdot x$ is an endomorphism. Let $\mathbf{B}/\theta$ be a finite subdirectly irreducible algebra in $\mathcal{V}(\mathbf{A})$ where **B** is a subalgebra of $\mathbf{A}^k$ and assume that this representation is selected so that $k$ is minimal. We claim that all $k$–tuples with at least one 0 coordinate are $\theta$–related to $\widehat{0}$. Suppose not, and let $\mathbf{b} = (b_1, \ldots, b_k) \in \mathbf{B}$ be such that $(\widehat{0}, \mathbf{b}) \notin \theta$ and **b** has a maximum number of 0 coordinates. We may assume without loss of generality that $b_1 = \ldots = b_l = 0$ ($l \geq 1$) and $b_{l+1}, \ldots, b_k \in A - \{0\}$. Consider the subuniverse $T = \{\mathbf{t} \in \mathbf{B}: t_1 = \ldots = t_l = 0\}$ of **B** and denote the corresponding subalgebra by **T**. The mapping

$$\rho \colon \mathbf{B} \to \mathbf{T}, \quad \mathbf{x} = (x_1, \ldots, x_k) \mapsto \mathbf{b} \cdot \mathbf{x} = (0, \ldots, 0, x_{l+1}, \ldots, x_k)$$

is an endomorphism of **B** onto **T**. Let $\alpha = \rho^{-1}(\theta|_T)$, let $\eta$ denote the kernel of the projection of **B** onto its last $k-1$ coordinates, and let $\tau$ be the equivalence relation on $B$ whose only nonsingleton block is $T$. It is straightforward to check that $\tau$ and $\beta = \theta \circ \tau \circ \theta$ are congruences of **B**. Clearly, $\theta \subseteq \alpha, \beta$, $(\widehat{0}, \mathbf{b}) \in \tau \subseteq \beta$, $(\widehat{0}, \mathbf{b}) \notin \theta$, $\eta \subseteq \alpha$, and by the minimality of $k$, $\eta \not\subseteq \theta$. Therefore $\theta \subset \alpha, \beta$. Since $\mathbf{B}/\theta$ is subdirectly irreducible, we conclude that $\theta \subset \alpha \cap \beta$. Let $(\mathbf{u}, \mathbf{v}) \in (\alpha \cap \beta) - \theta$. Then $\mathbf{b} \cdot \mathbf{u} \ \theta \ \mathbf{b} \cdot \mathbf{v}$ and $\mathbf{u} \ \theta \ \mathbf{s} \ \tau \ \mathbf{t} \ \theta \ \mathbf{v}$ for some $\mathbf{s}, \mathbf{t} \in T$. Hence

$$\mathbf{u} \ \theta \ \mathbf{s} = \mathbf{b} \cdot \mathbf{s} \ \theta \ \mathbf{b} \cdot \mathbf{u} \ \theta \ \mathbf{b} \cdot \mathbf{v} \ \theta \ \mathbf{b} \cdot \mathbf{t} = \mathbf{t} \ \theta \ \mathbf{v},$$

which is impossible as $(\mathbf{u}, \mathbf{v}) \notin \theta$. This contradiction finishes the proof of the claim that all elements of $\mathbf{B}$ having at least one 0 coordinate are $\theta$–related.

Since $\mathbf{B}/\theta$ is not a one–element algebra, it contains an element $\mathbf{c}/\theta$ with $\mathbf{c} \in (A - \{0\})^k$. Then $\widehat{2}/\theta = e(\mathbf{c})/\theta \in \mathbf{B}/\theta$ and $\widehat{1}/\theta = f(\widehat{2})/\theta \in \mathbf{B}/\theta$. It follows that $\widehat{0}/\theta, \widehat{1}/\theta, \widehat{2}/\theta$ are pairwise distinct. Thus $|\mathbf{B}/\theta| \geq 3$, the elements of $\mathbf{B}/\theta$ are

$$\widehat{0}/\theta, \ \widehat{1}/\theta, \ \widehat{2}/\theta, \ \mathbf{c}/\theta, \ \mathbf{d}/\theta, \ \ldots \quad \text{with} \quad \mathbf{c}, \mathbf{d}, \ldots \in (A - \{0\})^k,$$

and the operations act as follows:

$$x \cdot y = \begin{cases} \widehat{0}/\theta & \text{if } x = \widehat{0}/\theta \\ y & \text{otherwise} \end{cases},$$

and

$$e(x) = \begin{cases} \widehat{0}/\theta & \text{if } x = \widehat{0}/\theta \\ \widehat{2}/\theta & \text{otherwise} \end{cases}, \quad f(x) = \begin{cases} \widehat{0}/\theta & \text{if } x \neq \widehat{2}/\theta \\ \widehat{1}/\theta & \text{if } x = \widehat{2}/\theta \end{cases}.$$

The equivalence relation with a unique nonsingleton block $\{\widehat{1}/\theta, \mathbf{c}/\theta\}$ is a congruence on $\mathbf{B}/\theta$. A similar statement holds with $\mathbf{d}$ in place of $\mathbf{c}$, and this contradicts the subdirect irreducibility unless $B/\theta \subseteq \{\widehat{0}/\theta, \widehat{1}/\theta, \widehat{2}/\theta, \mathbf{c}/\theta\}$. From this and the description of the operations one deduces that $\mathbf{B} \cong \mathbf{A}$ or $\mathbf{B} \cong \mathbf{S}$.

As was mentioned before, $\mathbf{A}$ is a simple algebra of type **5**. One can show that $\mathbf{S}$ is subdirectly irreducible, its monolith $\mu$ is the equivalence relation collapsing $\widehat{1}/\gamma$ and $(1, 2)/\gamma$, and $\text{typ}(0, \mu) = \mathbf{1}$. In particular, the algebra $\mathbf{A}$ of this example omits types **1** and **2** and is self–rectangulating and yet it generates a variety which is not of type–set $\{\mathbf{5}\}$.

## 3  A Compatible Semilattice Term

If $\mathcal{V}$ is a locally finite variety, we have defined what it means for $\mathcal{V}$ to be self–rectangulating. When $\mathcal{V}$ has this property and $\mathcal{V}$ omits types **1** and **2**, then we must have $\text{typ}\{\mathcal{V}\} \subseteq \{\mathbf{5}\}$. The purpose of this section is to prove the following theorem.

**THEOREM 3.1** *If $\mathcal{V}$ is a locally finite variety of type set $\{\mathbf{5}\}$ then $\mathcal{V}$ is self–rectangulating if and only if it has a compatible semilattice term operation.*

We already know one half of this statement to be true: if $\mathcal{V}$ has a compatible semilattice term operation $\wedge$, then $\wedge$ interprets in every algebra as a compatible semilattice operation. The equivalence of the conditions of Theorem 2.8 completes the argument that $\mathcal{V}$ is self–rectangulating. In this section we shall concentrate on proving the reverse: that if $\mathcal{V}$ is self–rectangulating and typ$\{\mathcal{V}\} = \{\mathbf{5}\}$, then $\mathcal{V}$ has a compatible semilattice term operation.

First we reduce the claim to idempotent varieties. For a variety $\mathcal{V}$ we define a variety $\text{Id}(\mathcal{V})$ as follows. The language $\mathcal{L}$ of $\text{Id}(\mathcal{V})$ is specified by saying that there is one basic operation symbol corresponding to every term operation $t(\mathbf{x})$ of $\mathcal{V}$ for which
$$\mathcal{V} \models t(x, x, \ldots, x) = x.$$
Each algebra $\mathbf{A} \in \mathcal{V}$ gives rise to a model of $\mathcal{L}$ by interpreting the basic operation symbols just as one would interpret the term operation that gave rise to the basic operation symbol. The model of $\mathcal{L}$ which results is the idempotent reduct of $\mathbf{A}$ provided with a fixed indexing of its operations. Let $\mathcal{I}$ denote the class of all models of $\mathcal{L}$ which arise from members of $\mathcal{V}$ in this way. We define $\text{Id}(\mathcal{V})$ to be $\mathsf{HSP}(\mathcal{I})$. Clearly, $\text{Id}(\mathcal{V})$ is an idempotent variety. Moreover, for every equation satisfied in $\mathcal{V}$ which involves only idempotent term operations there is a corresponding equation satisfied in $\text{Id}(\mathcal{V})$.

**LEMMA 3.2** *If $\mathcal{V}$ is a locally finite variety of type–set $\{\mathbf{5}\}$ which is self–rectangulating, then so is $\text{Id}(\mathcal{V})$, and the variety $\mathcal{V}$ has a compatible semilattice term operation if and only if $\text{Id}(\mathcal{V})$ has a semilattice term operation.*

*Proof.* Assume that $\mathcal{V}$ is a locally finite variety of type–set $\{\mathbf{5}\}$ which is self–rectangulating. Since $\mathcal{V}$ is locally finite, $\mathcal{I}$ is uniformly locally finite and therefore $\text{Id}(\mathcal{V})$ is locally finite. Since $\mathcal{V}$ is of type–set $\{\mathbf{5}\}$, then by Theorem 9.10 of [2] $\mathcal{V}$ satisfies an idempotent Mal'cev condition which for locally finite varieties is equivalent to omitting types $\mathbf{1}$ and $\mathbf{2}$. This Mal'cev condition is expressible with equations involving idempotent term operations, so $\text{Id}(\mathcal{V})$ satisfies the same Mal'cev condition. Hence, $\text{Id}(\mathcal{V}) \cap \{\mathbf{1}, \mathbf{2}\} = \emptyset$. Finally, since $\mathcal{V}$ is self–rectangulating of type–set $\{\mathbf{5}\}$, we have by Theorem 2.8 that every finite member of $\mathcal{V}$ has an extension which has a compatible semilattice operation. This holds for the finite algebras in $\mathcal{I}$ and therefore for the finite algebras in $\mathsf{SP}(\mathcal{I})$. Hence Theorem 2.8 can be applied to prove that the finite free algebras of $\text{Id}(\mathcal{V})$ are self–rectangulating. It follows that all finite algebras in $\text{Id}(\mathcal{V})$ are self–rectangulating. Since $\text{Id}(\mathcal{V}) \cap \{\mathbf{1}, \mathbf{2}\} = \emptyset$ and $\text{Id}(\mathcal{V})$ is self–rectangulating, we get that typ$\{\text{Id}(\mathcal{V})\} = \{\mathbf{5}\}$. This proves the first claim of the lemma.

The clone of $\text{Id}(\mathcal{V})$ is isomorphic to the clone of idempotent term operations of $\mathcal{V}$. It is clear from this that if $\mathcal{V}$ has a compatible semilattice term operation then $\text{Id}(\mathcal{V})$ has one. Conversely, assume that $\text{Id}(\mathcal{V})$ has a semilattice term operation $t$ (not assumed to be compatible). Clearly, $t$ is a semilattice term operation of $\mathcal{V}$. By Theorem 2.8 every finite member $\mathbf{A}$ of $\mathcal{V}$ has a finite extension $\mathbf{B} \in \mathcal{V}$ which has a compatible semilattice operation $\wedge$. Now Lemma 2.5 shows that the term operation $t$ on $\mathbf{B}$ coincides with $\wedge$. Consequently $t$ is a compatible semilattice term operation of $\mathbf{B}$, and hence of $\mathbf{A}$. Since $\mathcal{V}$ is locally finite, this implies that $t$ is a compatible semilattice term operation of $\mathcal{V}$. □

From now on we can concentrate on idempotent varieties. Using a construction different from the one in the previous section we show that every self–rectangulating finite idempotent algebra $\mathbf{A}$ can be extended to a finite algebra $\mathbf{B} \in \mathcal{V}(\mathbf{A})$ such that $\mathbf{B}$ has a compatible semilattice operation and, in addition, the corresponding semilattice has a top element. As in the previous section, we start our considerations with subdirectly irreducible algebras.

**LEMMA 3.3** *Let $\mathbf{A}$ be a finite subdirectly irreducible algebra with non-abelian monolith $\mu$. Consider a $\langle 0_{\mathbf{A}}, \mu \rangle$–gene $G = (e, \sqcap, 0, 1, 0_{\mathbf{A}})$ where $U = e(A)$ is a $\langle 0_{\mathbf{A}}, \mu \rangle$–minimal set, and let $\sqsubseteq$ be the natural order of $\mathbf{A}$ with respect to $\langle 0_{\mathbf{A}}, \mu \rangle$. If $\mathbf{A}$ rectangulates itself with respect to $\langle 0_{\mathbf{A}}, \mu \rangle$, then*

(1) $x \wedge y = (x \sqcap y)|_U$ *is a compatible semilattice operation of the induced minimal algebra $\mathbf{A}|_U$;*

(2) $\sqsubseteq|_U$ *coincides with the semilattice order $\leq$ of $(U; \wedge)$;*

(3) *for every n-ary polynomial operation $p$ of $\mathbf{A}$ whose range is contained in $U$ and for all elements $x_1, \ldots, x_n, y_1, \ldots, y_n \in A$ we have*

$$p(x_1, \ldots, x_n) \wedge p(y_1, \ldots, y_n) = \bigwedge_{i=1}^{n} p(y_1, \ldots, y_{i-1}, x_i, y_{i+1}, \ldots, y_n).$$

*Proof.* First we observe that $\underline{1}$ is the largest element of $U$ with respect to the order $\sqsubseteq|_U$. For any element $u \in U$ we have $u \sqcap \underline{1} = u = \underline{1} \sqcap u$. Therefore, whenever $ef(u) = \underline{1}$ holds for some $f \in \mathrm{Pol}_1(\mathbf{A})$, then applying self–rectangulation to

$$ef(u \sqcap \underline{1}) = ef(u) = \underline{1} = ef(\underline{1} \sqcap u)$$

we obtain that $\underline{1} = ef(\underline{1} \sqcap \underline{1}) = ef(\underline{1})$. This shows that $u \sqsubseteq \underline{1}$.

To prove that the operation $x \wedge y = (x \sqcap y)|_U$ is commutative, consider arbitrary elements $u, v \in U$. By symmetry it suffices to show that $u \sqcap v \sqsubseteq v \sqcap u$. Assume that for some $f \in \mathrm{Pol}_1(\mathbf{A})$ we have $ef(u \sqcap v) = \underline{1}$. Then

$$\underline{1} = ef(u \sqcap v) \sqsubseteq ef(u \sqcap \underline{1}) = ef(u) = ef(u \sqcap u),$$

hence $\underline{1} = ef(u \sqcap u)$. Similarly, we get that $\underline{1} = ef(v \sqcap v)$. Thus

$$ef(u \sqcap \underline{u}) = \underline{1} = ef(\underline{v} \sqcap v),$$

which implies by self–rectangulation that $ef(v \sqcap u) = \underline{1}$, as required.

The associativity of $x \wedge y = (x \sqcap y)|_U$ can be verified in a similar fashion. Let $u, v, w \in U$. By commutativity it is enough to show that $u \sqcap (v \sqcap w) \sqsubseteq w \sqcap (v \sqcap u)$. Suppose that $ef(u \sqcap (v \sqcap w)) = \underline{1}$ for some $f \in \mathrm{Pol}_1(\mathbf{A})$. We have

$$\begin{aligned}\underline{1} = ef(u \sqcap (v \sqcap w)) &\sqsubseteq ef(u \sqcap (v \sqcap \underline{1})) \\ &= ef(u \sqcap v) \\ &= ef(u \sqcap (u \sqcap v)) \\ &= ef(u \sqcap (v \sqcap u))\end{aligned}$$

and
$$1 = ef(u \sqcap (v \sqcap w)) \sqsubseteq ef(1 \sqcap (v \sqcap w))$$
$$= ef(v \sqcap w)$$
$$= ef(w \sqcap v)$$
$$= ef(w \sqcap (w \sqcap v))$$
$$= ef(w \sqcap (v \sqcap w)).$$

Consequently,
$$ef(u \sqcap \underline{(v \sqcap u)}) = 1 = ef(\underline{w} \sqcap (v \sqcap w)),$$
implying by self–rectangulation that $ef(w \sqcap (v \sqcap u)) = 1$.

Thus $x \wedge y = (x \sqcap y)|_U$ is a semilattice operation on $U$. Obviously, for any elements $u, v \in U$ we have $u \wedge v \sqsubseteq u \wedge 1 = u$ and similarly $u \wedge v \sqsubseteq v$. Furthermore, if $w \sqsubseteq u, v$ ($w \in U$), then $w = w \wedge w \sqsubseteq u \wedge v$, showing that $u \wedge v$ is the greatest lower bound of $u$ and $v$ with respect to $\sqsubseteq|_U$. This proves (2).

To conclude the proof of (1) it remains to show that $x \wedge y$ is compatible with the operations of $\mathbf{A}|_U$. Let $p$ be an $n$-ary polynomial operation of $\mathbf{A}$ whose range is contained in $U$, and let $\mathbf{u}, \mathbf{v} \in U^n$. We claim that

$$p(\mathbf{u}) \wedge p(\mathbf{v}) = p(\mathbf{u} \wedge \mathbf{v}).$$

By (2) the inequality $\sqsupseteq$ is obvious. To prove the reverse inequality assume that $ef(p(\mathbf{u}) \wedge p(\mathbf{v})) = 1$ for some $f \in \mathrm{Pol}_1(\mathbf{A})$. It follows that $ef(p(\mathbf{u})) \sqsupseteq 1$ and $ef(p(\mathbf{v})) \sqsupseteq 1$. Hence

$$ef(p(\underline{\mathbf{u}} \wedge \mathbf{u})) = ef(p(\mathbf{u})) = 1 = ef(p(\mathbf{v})) = ef(p(\mathbf{v} \wedge \underline{\mathbf{v}})),$$

so by self–rectangulation we conclude that $ef(p(\mathbf{u} \wedge \mathbf{v})) = 1$.

The equality in (3) can be established by showing that for any $f \in \mathrm{Pol}_1(\mathbf{A})$ the $ef$–image of the left hand side equals 1 if and only if the $ef$–image of the right hand side equals 1. We know from (1)–(2) that $ef|_U$ is a $\wedge$–endomorphism and 1 is the largest element of the semilattice $(U; \wedge)$. Therefore the $ef$–image of the left hand side of the equality is 1 if and only if

$$ef(p(x_1, \ldots, x_n)) = 1 = ef(p(y_1, \ldots, y_n)),$$

and the $ef$–image of the right hand side of the equality is 1 if and only if

$$ef(p(y_1, \ldots, y_{i-1}, x_i, y_{i+1}, \ldots, x_n)) = 1 \quad \text{for all } 1 \leq i \leq n.$$

Making use of self–rectangulation one can easily see that the latter two conditions are equivalent. $\square$

**LEMMA 3.4** *Let $\mathbf{A}$ be a finite subdirectly irreducible algebra with nonabelian monolith $\mu$ such that $\mathbf{A}$ rectangulates itself with respect to $\langle 0_{\mathbf{A}}, \mu \rangle$, and let $U$ be a $\langle 0_{\mathbf{A}}, \mu \rangle$–minimal set. Furthermore, let $\wedge$ be the (unique) compatible semilattice term operation of $\mathbf{A}|_U$ (cf. Lemma 3.3 and Corollary 2.6). If the basic operations of the algebra $\mathbf{A}$ are surjective, then $\mathbf{A}$ is isomorphic to a subreduct of a matrix power of the semilattice $(U; \wedge)$.*

*Proof.* Let $G = (e, \sqcap, 0, 1, 0_{\mathbf{A}})$ be a $\langle 0_{\mathbf{A}}, \mu \rangle$–gene with $U = e(A)$ as in Lemma 3.3. We define a set $F$ of unary polynomial operations of $\mathbf{A}$ as follows:

$$F = \{ef \colon f \in \mathrm{Pol}_1(\mathbf{A}),\ 1 \in ef(A)\},$$

and enumerate $F$ as $F = \{e_1, \ldots, e_m\}$.

Now we consider the mapping

$$\varphi \colon A \to U^m, \quad a \mapsto (e_1(a), \ldots, e_m(a)).$$

If for some elements $a, b \in A$ we have $e_i(a) = e_i(b)$ for all $i$ ($1 \leq i \leq m$), then the definition of the natural order $\sqsubseteq$ immediately implies that $a = b$. Thus $\varphi$ is injective.

Let $g(x_1, \ldots, x_n)$ be any basic operation of $\mathbf{A}$. For any $i$ ($1 \leq i \leq m$) the range of $e_i g(x_1, \ldots, x_n)$ has a largest element, namely 1, because $g$ is surjective and $1 \in e_i(A) \subseteq U$. Select elements $a_{i1}, \ldots, a_{in} \in A$ such that $e_i g(a_{i1}, \ldots, a_{in}) = 1$. Applying Lemma 3.3 we see that

$$\begin{aligned} e_i g(x_1, \ldots, x_n) &= e_i g(x_1, \ldots, x_n) \wedge e_i g(a_{i1}, \ldots, a_{in}) \\ &= \bigwedge_{j=1}^{n} e_i g(a_{i1}, \ldots, a_{i,j-1}, x_j, a_{i,j+1}, \ldots, a_{in}) \end{aligned}$$

for all $x_1, \ldots, x_n \in A$. The meetands on the right hand side are unary polynomial operations of $\mathbf{A}$ whose range is contained in $U$; moreover, this range contains 1 because $e_i g(a_{i1}, \ldots, a_{in}) = 1$. Thus the meetands are members of $F$. In symbols, we have

$$e_i g(x_1, \ldots, x_n) = \bigwedge_{j=1}^{n} e_{\lambda(i,j)}(x_j)$$

for appropriate indices $\lambda(i, j) \in \{1, \ldots, m\}$.

It is straightforward to check that if we make correspond to each $n$–ary basic operation $g$ of $\mathbf{A}$ the $n$-ary operation

$$(y_1^1, \ldots, y_1^m), \ldots, (y_n^1, \ldots, y_n^m) \overset{\bar{g}}{\mapsto} \left( \bigwedge_{j=1}^{n} y_j^{\lambda(1,j)}, \ldots, \bigwedge_{j=1}^{n} y_j^{\lambda(m,j)} \right)$$

on $U^m$, then $\varphi$ yields an isomorphism between $\mathbf{A} = (A; g, \ldots)$ and a subalgebra of the reduct $\mathbf{R} = (U^m; \bar{g}, \ldots)$ of $(U; \wedge)^{[m]}$. □

**Added later:** Since this was written R. McKenzie discovered a characterization of those finite algebras which are isomorphic to a subalgebra of a reduct of a matrix power of the two-element lower bounded meet semilattice. His result can be applied to prove Lemma 3.4 with the surjectivity assumption omitted.

**THEOREM 3.5** *Let $\mathbf{A}$ be a finite idempotent algebra which omits types $\mathbf{1}$ and $\mathbf{2}$. If $\mathbf{A}$ is self–rectangulating, then $\mathbf{A}$ is isomorphic to a subalgebra of a finite algebra $\mathbf{B} \in \mathcal{V}(\mathbf{A})$ such that $\mathbf{B}$ has a compatible semilattice operation and the corresponding semilattice has a largest element.*

*Proof.* Every finite algebra $\mathbf{A}$ which omits types **1** and **2** and is self–rectangulating is a subdirect product of subdirectly irreducible algebras having all of these properties. It suffices to prove the claim of the theorem for subdirectly irreducible algebras, since the class of algebras $\mathbf{A}$ which have extensions in $\mathcal{V}(\mathbf{A})$ with the required properties is closed under the formation of subalgebras and products.

From now on we assume that $\mathbf{A}$ is subdirectly irreducible, and we use the notation introduced in the statement and the proof of the preceding lemma. Let $\mathbf{W}$ denote the subalgebra of the reduct $\mathbf{R}$ of $(U; \wedge)^{[m]}$ with base set $\varphi(A)$ which is isomorphic to $\mathbf{A}$ via $\varphi$. Define a subset $B$ of $U^m$ as follows:

$$B = \Big\{(y^1, \ldots, y^m) \in U^m \colon \text{ for all } L \subseteq \{1, \ldots, m\} \text{ and } 1 \leq r \leq m,$$
$$\bigwedge_{l \in L} y^l = y^r \text{ whenever } \mathbf{A} \models \bigwedge_{l \in L} e_l(x) = e_r(x)\Big\}.$$

We show that $B$ is a subuniverse of $\mathbf{R}$. Let $g$ be a basic operation of $\mathbf{A}$, and use all notation related to $g$ which was introduced in the previous proof; in particular, $\bar{g}$ is the corresponding operation of $\mathbf{R}$. Select arbitrary $m$-tuples $y_j = (y_j^1, \ldots, y_j^m) \in B$ ($1 \leq j \leq n$). To verify that $\bar{g}(y_1, \ldots, y_n) \in B$, we consider an arbitrary equation

$$\mathbf{A} \models \bigwedge_{l \in L} e_l(x) = e_r(x)$$

where $L \subseteq \{1, \ldots, m\}$ and $1 \leq r \leq m$. Then

$$\mathbf{A} \models \bigwedge_{l \in L} e_l g(x_1, \ldots, x_n) = e_r g(x_1, \ldots, x_n),$$

that is

$$\mathbf{A} \models \bigwedge_{l \in L} \bigwedge_{j=1}^{n} e_{\lambda(l,j)}(x_j) = \bigwedge_{j=1}^{n} e_{\lambda(r,j)}(x_j). \tag{1}$$

Since $g$ is idempotent, we have that

$$\mathbf{A} \models \bigwedge_{j=1}^{n} e_{\lambda(r,j)}(x) = e_r g(x, \ldots, x) = e_r(x),$$

whence it follows that $e_{\lambda(l,j)}(x) \geq e_r(x)$ and $e_{\lambda(r,j)}(x) \geq e_r(x)$ for all $l \in L$, $1 \leq j \leq n$, and $x \in A$. By construction, $e_r(a) = 1$ for some $a \in A$. Since 1 is the largest element of $(U; \wedge)$, this implies that $e_{\lambda(l,j)}(a) = e_{\lambda(r,j)}(a) = 1$ for all $l \in L$ and $1 \leq j \leq n$. Consequently, if we substitute $a$ for all variables in (1) but one, then we conclude that

$$\mathbf{A} \models \bigwedge_{l \in L} e_{\lambda(l,j)}(x) = e_{\lambda(r,j)}(x) \qquad \text{for all } 1 \leq j \leq n.$$

This shows that in (1) equality holds separately in each variable. Hence by the definition of $B$ we have

$$\bigwedge_{l \in L} y_j^{\lambda(l,j)} = y_j^{\lambda(r,j)} \qquad \text{for all } 1 \leq j \leq n.$$

By taking $\bigwedge$ over all $j$ we see that the required equality holds for the coordinates of $\bar{g}(y_1, \ldots, y_n)$. This proves that $\bar{g}(y_1, \ldots, y_n) \in B$.

Let $\mathbf{B}$ denote the subalgebra of $\mathbf{R}$ with base set $B$. Clearly, $\mathbf{W}$ ($\cong \mathbf{A}$) is a subalgebra of $\mathbf{B}$. Furthermore, $B$ is closed under the coordinatewise action of the operation $x \wedge y$, and this yields a compatible semilattice operation for $\mathbf{B}$. The $m$-tuple $(1, \ldots, 1)$ belongs to $\mathbf{B}$, and it is a largest element with respect to this semilattice operation, as $B \subseteq U^m$ and 1 is the largest element of the semilattice $(U; \wedge)$.

It remains to show that $\mathbf{B} \in \mathcal{V}(\mathbf{A})$. Let $t, t'$ be $n$-ary terms in the language of $\mathbf{A}$ such that $\mathbf{A} \models t = t'$, that is, $\mathbf{W} \models t = t'$. By construction $\mathbf{W}$ is a subalgebra of $\mathbf{R}$, and $\mathbf{R}$ is a reduct of $(U; \wedge)^{[m]}$. Therefore the term operation corresponding to $t$ has the form

$$(y_1^1, \ldots, y_1^m), \ldots, (y_n^1, \ldots, y_n^m) \mapsto \left( \ldots, \overbrace{\bigwedge_{(k,j) \in T_i} y_j^k}^{i\text{th}}, \ldots \right) \tag{2}$$

where $T_i$ ($1 \leq i \leq m$) are appropriate subsets of $\{1, \ldots, m\} \times \{1, \ldots, n\}$, and the term operation corresponding to $t'$ has a similar form with $T_i'$ in place of $T_i$ for each $i$.

The way the isomorphism between $\mathbf{A}$ and $\mathbf{W}$ was set up shows that for $1 \leq i \leq m$ we have

$$\mathbf{A} \models e_i t(x_1, \ldots, x_n) = \bigwedge_{(k,j) \in T_i} e_k(x_j), \quad e_i t'(x_1, \ldots, x_n) = \bigwedge_{(k,j) \in T_i'} e_k(x_j).$$

Hence the assumption $\mathbf{A} \models t = t'$ is equivalent to the condition that

$$\mathbf{A} \models \bigwedge_{(k,j) \in T_i} e_k(x_j) = \bigwedge_{(k,j) \in T_i'} e_k(x_j) \quad \text{for all } 1 \leq i \leq m.$$

The same argument as before implies that equality holds separately in each variable; that is, we have

$$\mathbf{A} \models \bigwedge_{\substack{k \\ (k,j) \in T_i}} e_k(x) = \bigwedge_{\substack{k \\ (k,j) \in T_i'}} e_k(x) \quad \text{for all } 1 \leq i \leq m, \ 1 \leq j \leq n.$$

Furthermore, since $\mathbf{A}$ is idempotent,

$$\bigwedge_{\substack{k \\ (k,j) \in T_i}} e_k(x) \geq e_i t(x, \ldots, x) = e_i(x) \quad \text{for all } x \in A, \ 1 \leq i \leq m, \ 1 \leq j \leq n.$$

By construction, each unary polynomial operation $e_i$ assumes the value 1; let $e_i(a_i) = 1$ ($a_i \in A$). Hence

$$\bigwedge_{\substack{k \\ (k,j) \in T_i}} e_k(x) = e_i t(a_1, \ldots, a_{j-1}, x, a_{j+1}, \ldots, a_n)$$

is a unary polynomial operation, and it also assumes the value 1. Therefore for all indices $1 \leq i \leq m$, $1 \leq j \leq n$ there exists an index $1 \leq l \leq m$ such that

$$\mathbf{A} \models \bigwedge_{\substack{k \\ (k,j) \in T_i}} e_k(x) = e_l(x) = \bigwedge_{\substack{k \\ (k,j) \in T'_i}} e_k(x).$$

By the definition of $B$ this implies that for all $m$-tuples $y = (y^1, \ldots, y^m) \in B$ and for all indices $1 \leq i \leq m$, $1 \leq j \leq n$ the equality

$$\bigwedge_{\substack{k \\ (k,j) \in T_i}} y^k = y^l = \bigwedge_{\substack{k \\ (k,j) \in T'_i}} y^k.$$

holds. If we form the meet of the left and right sides over all $j$ we obtain that $\bigwedge_{(k,j) \in T_i} y_j^k = \bigwedge_{(k,j) \in T'_i} y_j^k$. Referring to (2) we see that this equality for each $i$ implies that $\mathbf{B} \models t = t'$. $\square$

**LEMMA 3.6** *Let $\mathbf{B}$ be a finite idempotent algebra such that $\mathbf{B}$ has a compatible semilattice operation $x \wedge y$ and the corresponding semilattice has a largest element. If $\mathbf{B}$ contains no 2-element essentially unary subalgebra, then $x \wedge y$ is a term operation of $\mathbf{B}$.*

*Proof.* The natural order and the largest element of the semilattice $(B; \wedge)$ will be denoted by $\leq$ and $1$, respectively. For $x, y \in \mathbf{B}$ and any binary term operation $t$ we have

$$\begin{aligned} t(x, y) &= t(x \wedge 1, 1 \wedge y) = t(x, 1) \wedge t(1, y), \\ t(x \wedge y, 1) &= t(x \wedge y, 1 \wedge 1) = t(x, 1) \wedge t(y, 1). \end{aligned}$$

Thus, $t$ has the form $t(x, y) = \alpha(x) \wedge \beta(y)$ for some unary polynomial operations $\alpha, \beta$ of $\mathbf{B}$, which are $\wedge$-endomorphisms. The idempotence of $\mathbf{B}$ implies that for all elements $x \in B$ we have $\alpha(x) \wedge \beta(x) = x$, and hence $\alpha(x), \beta(x) \geq x$. This shows that any binary term operation is expressible as $\alpha(x) \wedge \beta(y)$ where $\alpha$ and $\beta$ are increasing $\wedge$–endomorphisms.

Fixing a binary term $t$ and a representation $t(x, y) = \alpha(x) \wedge \beta(y)$, we have

$$\alpha(\alpha^n(x) \wedge \beta(y)) \wedge \beta(y) = \alpha^{n+1}(x) \wedge \alpha\beta(y) \wedge \beta(y) = \alpha^{n+1}(x) \wedge \beta(y)$$

for any $n$, since $\alpha$ is an increasing $\wedge$-endomorphism. This implies that if we iterate $t(x, y)$ in its first variable until we produce a term $t'(x, y)$ which satisfies $t'(t'(x, y), y) = t'(x, y)$, then the term operation $t'(x, y)$ may be represented as $\bar{\alpha}(x) \wedge \beta(y)$ where $\bar{\alpha}^2 = \bar{\alpha}$. We can repeat this argument in the second variable to produce a term $t''(x, y)$ whose representation is $\bar{\alpha}(x) \wedge \bar{\beta}(y)$ where $\bar{\alpha}^2 = \bar{\alpha}$ and $\bar{\beta}^2 = \bar{\beta}$. If $s(x, y) = t''(t''(x, y), t''(y, x))$, then $s(x, y)$ has the representation

$$\begin{aligned} s(x, y) &= \bar{\alpha}(\bar{\alpha}(x) \wedge \bar{\beta}(y)) \wedge \bar{\beta}(\bar{\alpha}(y) \wedge \bar{\beta}(x)) \\ &= (\bar{\alpha}(x) \wedge \bar{\beta}(x)) \wedge (\bar{\alpha}\bar{\beta}(y) \wedge \bar{\beta}\bar{\alpha}(y)) \\ &= t''(x, x) \wedge (\bar{\alpha}\bar{\beta}(y) \wedge \bar{\beta}\bar{\alpha}(y)) \\ &= x \wedge (\bar{\alpha}\bar{\beta}(y) \wedge \bar{\beta}\bar{\alpha}(y)). \end{aligned}$$

That is, $s(x,y)$ is a term operation of **B** of the form $x \wedge \gamma(y)$ where $\gamma$ is an increasing $\wedge$-endomorphism. Iteration in the second variable yields a term operation $x \wedge \varepsilon(y)$ such that $\varepsilon^2 = \varepsilon$.

Let us enumerate the term operations of **B** of the form $x \wedge \varepsilon(y)$ with $\varepsilon^2 = \varepsilon$ as $x \wedge \varepsilon_i(y)$ ($\varepsilon_i^2 = \varepsilon_i$), $i = 1, 2, \ldots, k$, and form the term operation

$$x \wedge \delta(y) = (\ldots ((x \wedge \varepsilon_1(y)) \wedge \varepsilon_2(y)) \wedge \ldots) \wedge \varepsilon_k(y).$$

As before, for some power $E = \delta^m$ of $\delta$ we have that $E^2 = E$. The $m$-th iterate of $x \wedge \delta(y)$ in the second variable is the term operation $x \wedge E(y)$. Clearly,

$$E(x) = \delta^m(x) \leq \varepsilon_i^m(x) = \varepsilon_i(x) \qquad \text{for all } x \in B,\ 1 \leq i \leq k.$$

We prove the contrapositive of the claim of the lemma. Assume $x \wedge y$ is not a term operation of **B**. Then $E$ is not the identity endomorphism of the semilattice $(B; \wedge)$. Hence there exists an element $u \in B$ such that $E(u) \neq u$. Let $v = E(u)$. Since $E(u) \geq u$, it follows that $v > u$. Furthermore, $E(v) = E^2(u) = E(u) = v$. This equality implies that every element in the interval $[u, v]$ is mapped by $E$ into $v$. Replacing $u$ by a lower cover of $v$ we can assume that $u$ is covered by $v$, and all properties established so far remain valid. In addition,

$$\varepsilon_i(u) \leq \varepsilon_i(v) = \varepsilon_i E(u) \leq \varepsilon_i^2(u) = \varepsilon_i(u),$$

therefore $\varepsilon_i(u) = \varepsilon_i(v)$ for all $i$ ($1 \leq i \leq k$). Since $u < v \leq \varepsilon_i(v)$, this equality shows that every term operation $x \wedge \varepsilon_i(y)$ ($1 \leq i \leq k$) restricts to $\{u, v\}$ as a projection.

**B** is idempotent, therefore for any binary term operation $t(x,y)$ of **B** we have $u \leq t(u,v) \leq v$. But $u$ is covered by $v$, hence $\{u,v\}$ is a subuniverse of **B**. Let **U** denote the subalgebra of **B** with base set $\{u,v\}$. The set $\{u,v\}$ is closed under the operation $x \wedge y$ as well, hence $(x \wedge y)|_{\{u,v\}}$ is a compatible semilattice operation of **U**. Thus **U** can have no other binary term operations than the projections and $(x \wedge y)|_{\{u,v\}}$. Suppose that $(x \wedge y)|_{\{u,v\}}$ is a term operation of **U**, say, $t(x,y)|_{\{u,v\}} = (x \wedge y)|_{\{u,v\}}$ for some binary term $t$. Going over the construction above yielding $x \wedge \varepsilon(y)$ from $t(x,y) = \alpha(x) \wedge \beta(y)$, one can easily check that in this case $(x \wedge \varepsilon(y))|_{\{u,v\}} = (x \wedge y)|_{\{u,v\}}$. However, this possibility was excluded in the preceding paragraph. Thus we have established that **U** is a 2–element algebra with a compatible semilattice operation such that every binary term operation of **U** is a projection. The only 2–element algebras with these properties are essentially unary. $\square$

Now we are ready to complete the proof of Theorem 3.1. Let $\mathcal{V}$ be a locally finite variety of type set $\{\mathbf{5}\}$ which is self–rectangulating. By Lemma 3.2 it will follow that $\mathcal{V}$ has a compatible semilattice term operation if we show that $\mathrm{Id}(\mathcal{V})$ has a semilattice term operation.

Let **F** be the two-generated free algebra in $\mathrm{Id}(\mathcal{V})$. By Lemma 3.2 **F** is a finite idempotent algebra of type set $\{\mathbf{5}\}$ which is self–rectangulating, therefore by Theorem 3.5 **F** can be extended to a finite algebra $\mathbf{B} \in \mathrm{Id}(\mathcal{V})$ such that **B** has a compatible semilattice operation and the corresponding semilattice has a largest element. Since $\mathrm{Id}(\mathcal{V})$ omits type $\mathbf{1}$, Lemma 3.6 implies that $\mathbf{B}$ has a compatible semilattice term operation. The same conclusion holds then for $\mathbf{F}$, too. Thus the variety $\mathrm{Id}(\mathcal{V})$ has a semilattice term operation, as required. $\square$

As a consequence of Theorem 3.1 we get the following improvement on Corollary 2.9.

**COROLLARY 3.7** *If $\mathbf{A}$ is a finite algebra which omits types $\mathbf{1}$ and $\mathbf{2}$, then the following conditions are equivalent.*

(i) $\mathbf{A}$ *is self–rectangulating and* $\mathcal{V}(\mathbf{A})$ *omits type* $\mathbf{1}$.

(ii) $\mathbf{A}$ *is self–rectangulating and* $\mathrm{typ}\{\mathcal{V}(\mathbf{A})\} = \{\mathbf{5}\}$.

(iii) $\mathbf{A}$ *has a compatible semilattice term operation.*

*Proof.* The statement that a locally finite variety $\mathcal{V}$ has a semilattice term implies that $\mathrm{typ}\{\mathcal{V}\} \cap \{\mathbf{1}, \mathbf{2}\} = \emptyset$ since this is true of the variety of semilattices. It is easily deducible from Theorem 2.8 that if $\mathbf{A}$ (and hence $\mathcal{V}(\mathbf{A})$) has a compatible semilattice term, then $\mathcal{V}(\mathbf{A})$ is self–rectangulating. Hence (iii)$\Rightarrow$(ii)$\Rightarrow$(i). The implication (i)$\Rightarrow$(iii) can be deduced from Corollary 2.9 and from Theorem 3.1. $\square$

## 4 Definable Principal Congruences

**LEMMA 4.1** *Let $\mathbf{A}$ be a finite algebra with a compatible semilattice term operation $x \wedge y$. If $t(x_1, \ldots, x_m)$ is a term in the language of $\mathbf{A}$, then for some $k \leq |A|^{|A|}$ there is a term $t'(y_1, \ldots, y_k)$ and a partition $\{X_1, \ldots, X_k\}$ of the set $\{x_1, \ldots, x_m\}$ such that*

$$\mathbf{A} \models t(x_1, \ldots, x_m) = t'(\bigwedge X_1, \ldots, \bigwedge X_k).$$

*Proof.* Let us fix a term $t$ and denote by $R$ the range of the term operation $t^{\mathbf{A}}$. If $\mathbf{c}, \mathbf{d} \in A^m$ are $m$–tuples such that $t^{\mathbf{A}}(\mathbf{c}) = t^{\mathbf{A}}(\mathbf{d})$ then $t^{\mathbf{A}}(\mathbf{c}) = t^{\mathbf{A}}(\mathbf{c}) \wedge t^{\mathbf{A}}(\mathbf{d}) = t^{\mathbf{A}}(\mathbf{c} \wedge \mathbf{d})$. Therefore for every element $a \in R$ there exists a smallest $m$-tuple $\mathbf{b}$ with $t^{\mathbf{A}}(\mathbf{b}) = a$; this $\mathbf{b}$ is the meet of all $\mathbf{c} \in A^m$ with $t^{\mathbf{A}}(\mathbf{c}) = a$.

Now we define an equivalence relation $\sim$ on $\{x_1, \ldots, x_m\}$ as follows: $x_i \sim x_j$ if and only if for every element $a \in R$, in the smallest $m$–tuple $\mathbf{b} = (b_1, \ldots, b_m)$ with $t^{\mathbf{A}}(\mathbf{b}) = a$ we have $b_i = b_j$. Clearly, $x_i \sim x_j$ holds exactly when the mappings $R \to A, a \mapsto b_i$ and $R \to A, a \mapsto b_j$ coincide. Therefore $\sim$ has at most $|A|^{|A|}$ blocks.

We show that if $x_1 \sim x_2$ then

$$\mathbf{A} \models t(x_1, \ldots, x_m) = t(x_1 \wedge x_2, x_1 \wedge x_2, x_3, \ldots, x_m);$$

that is, for the $(m-1)$-ary term $t'(y_1,\ldots,y_{m-1}) = t(y_1,y_1,\ldots,y_{m-1})$ we have

$$\mathbf{A} \models t(x_1,\ldots,x_m) = t'(x_1 \wedge x_2, x_3 \ldots, x_m).$$

Choose an arbitrary $m$-tuple $(c_1,\ldots,c_m) \in A^m$, and let $t^{\mathbf{A}}(c_1,\ldots,c_m) = a$. For the smallest $m$–tuple $(b_1,\ldots,b_m)$ with $t^{\mathbf{A}}(b_1,\ldots,b_m) = a$ we have $b_1 = b_2$. Therefore $b_1 = b_2 \leq c_1, c_2$, whence $b_1 = b_2 \leq c_1 \wedge c_2$. Thus

$$a = t^{\mathbf{A}}(b_1,\ldots,b_m) \leq t^{\mathbf{A}}(c_1 \wedge c_2, c_1 \wedge c_2, c_3 \ldots, c_m) \leq t^{\mathbf{A}}(c_1,\ldots,c_m) = a,$$

which implies the necessary equality $t^{\mathbf{A}}(c_1 \wedge c_2, c_1 \wedge c_2, c_3 \ldots, c_m) = a$.

The estimate on the number of blocks of $\sim$ and the claim proved in the preceding paragraph show that for any term $t$ of arity $m > |A|^{|A|}$ there exists a term $t'$ of smaller arity $k$ for which an identity required in the lemma holds. Applying this fact to $t'$ in place of $t$ we get that the least $k$ for which such a term $t'$ exists for $t$ must satisfy $k \leq |A|^{|A|}$. This completes the proof. □

**COROLLARY 4.2** *Assume that $\mathbf{A}$ is a finite algebra with a compatible semilattice term operation and that $\mathbf{B} \in \mathcal{V}(\mathbf{A})$. If $p(x)$ is a unary polynomial of $\mathbf{B}$, then for some $k \leq |A|^{|A|}$ there is an $(k+1)$–ary term $r$ and a tuple $\mathbf{b} \in B^k$ such that $p(x) = r^{\mathbf{B}}(x, \mathbf{b})$.*

*Proof.* Since $p \in \mathrm{Pol}_1(\mathbf{B})$ there is a term $t(x, \mathbf{y})$ and a tuple $\mathbf{c}$ of elements of $B$ such that $p(x) = t^{\mathbf{B}}(x, \mathbf{c})$. Using Lemma 4.1, there is $k \leq |A|^{|A|}$, a term $t'(x_1,\ldots,x_k)$ and a partition of the set $\{x, c_1,\ldots,c_m\}$ into $k$ subsets $X_1,\ldots,X_k$ such that

$$(t')^{\mathbf{B}}(\bigwedge X_1, \ldots, \bigwedge X_k) = t^{\mathbf{B}}(x, \mathbf{c}) = p(x).$$

We may assume that $X_1 = \{x, c_{11},\ldots,c_{1s_1}\}$ and that $X_i = \{c_{i1},\ldots,c_{is_i}\}$. Set $r(x, y_1,\ldots,y_k) = t'(x \wedge y_1, y_2, \ldots, y_k)$ and set $b_i = \bigwedge_j c_{ij}$. Then we have $p(x) = r^{\mathbf{B}}(x, b_1, \ldots, b_k)$ as desired. □

Let $\mathcal{L}$ be a language of algebras. A **principal congruence formula** of $\mathcal{L}$ is a positive existential first–order formula $\Phi(x, y, u, v)$ of $\mathcal{L}$ for which

$$\models \Phi(x, y, u, u) \Longrightarrow x = y$$

holds. If $\mathbf{A}$ is a model of $\mathcal{L}$ and $\Phi$ is a principal congruence formula of $\mathcal{L}$, then for $a, b, c, d \in A$ we have that $\Phi^{\mathbf{A}}(a, b, c, d)$ implies that $(a, b) \in \mathrm{Cg}^{\mathbf{A}}(c, d)$. (To see this, use the displayed implication in $\mathbf{A}/\mathrm{Cg}^{\mathbf{A}}(c, d)$.) We will call a formula a **special principal congruence formula** if it has the following form:

$$\exists \mathbf{a}, \mathbf{w}_1,\ldots,\mathbf{w}_{n-1}\, (x = a_1) \& (y = a_n) \&_{i=1}^{n-1} (p_i(\{u, v\}, \mathbf{w}_i) = \{a_i, a_{i+1}\}).$$

In this formula, $\mathbf{a} = (a_1,\ldots,a_n)$, the $p_i$ are terms in the language, and the expression $p_i(\{u, v\}, \mathbf{w}_i) = \{a_i, a_{i+1}\}$ means that $(p_i(u, \mathbf{w}_i) = a_i \& p_i(v, \mathbf{w}_i) = a_{i+1})$ or that the same thing holds with $u$ and $v$ switched.

It is easy to see that a special principal congruence formula is a principal congruence formula. Conversely, Mal'cev's congruence generation theorem proves that if $\mathbf{A}$ is an algebra with $a, b, c, d \in A$ and $(a, b) \in \mathrm{Cg}^{\mathbf{A}}(c, d)$, then there is a special principal congruence formula $\Phi$ such that $\Phi^{\mathbf{A}}(a, b, c, d)$ holds. It follows that for models of $\mathcal{L}$ the statement that $(x, y) \in \mathrm{Cg}(u, v)$ is equivalent to the (infinite) disjunction of all special principal congruence formulas of $\mathcal{L}$. If the 4–ary relation $(x, y) \in \mathrm{Cg}(u, v)$ is uniformly definable throughout a variety $\mathcal{V}$, then the compactness theorem implies that there is a principal congruence formula $\Phi$, which may be taken to be a finite disjunction of special principal congruence formulas, such that

$$\mathbf{A} \models (a, b) \in \mathrm{Cg}^{\mathbf{A}}(c, d) \iff \Phi^{\mathbf{A}}(a, b, c, d)$$

holds for all $\mathbf{A} \in \mathcal{V}$. When this happens we say that $\mathcal{V}$ has **definable principal congruences** or **DPC**. DPC is a rare but useful property for a variety. We refer the reader to [1] for an early investigation of DPC.

**THEOREM 4.3** *If $\mathbf{A}$ is a finite algebra which has a compatible semilattice term operation, then $\mathbf{A}$ generates a variety with DPC.*

*Proof.* Choose an algebra $\mathbf{B} \in \mathcal{V} = \mathcal{V}(\mathbf{A})$ and elements $a, b, c, d \in B$ such that $(a, b) \in \mathrm{Cg}(c, d)$. We shall argue that there is a special principal congruence formula $\Phi$ of a very restricted form for which $\Phi^{\mathbf{B}}(a, b, c, d)$ holds.

Choose any Mal'cev chain $a = a_1, \ldots, a_n = b$ which connects $a$ to $b$ by a chain of polynomial images of $c$ and $d$. Say that the polynomials involved are $p_i^{\mathbf{B}}(x, \mathbf{w}_i)$, $i = 1, \ldots, n$, and that $p_i^{\mathbf{B}}(\{c, d\}, \mathbf{w}_i) = \{a_i, a_{i+1}\}$. We alter this to a new chain by defining $g_i = a_1 \wedge \cdots \wedge a_i$ and $h_i = a_i \wedge \cdots \wedge a_n$. This gives us a descending chain $a = g_1, \ldots, g_n$ followed by an ascending chain $g_n = h_1, \ldots, h_n = b$. The polynomials which witness that this is indeed a Mal'cev chain are defined as follows: $r_i(x) = p_i^{\mathbf{B}}(x, \mathbf{w}_i) \wedge g_i$ and $s_i(x) = p_i^{\mathbf{B}}(x, \mathbf{w}_i) \wedge h_{i+1}$. Then we have $r_i(\{c, d\}) = \{g_i, g_{i+1}\}$ and $s_i(\{c, d\}) = \{h_i, h_{i+1}\}$.

Let $P$ be a minimal set of $\mathcal{V}$–inequivalent representatives for the $(|A|^{|A|}+1)$–ary terms. We will prefer the first variable in any member of $P$, so a typical member of $P$ will be written $t(x, \mathbf{y})$. From Corollary 4.2 we have that each polynomial $r_i(x)$ is equal to one of the form $R_i^{\mathbf{B}}(x, \mathbf{p}_i)$ where $R_i \in P$ and $\mathbf{p}_i \in B^{|A|^{|A|}}$. Similarly, each $s_i(x)$ is equal to one of the form $S_i^{\mathbf{B}}(x, \mathbf{q}_i)$. Furthermore, we shall argue, no two nontrivial links in the descending chain $g_1, \ldots, g_n$ can arise from the same term $R(x, \mathbf{y})$. To see this, assume that $g_i > g_{i+1} \geq g_j > g_{j+1}$ and that $R^{\mathbf{B}}(\{c, d\}, \mathbf{p}_i) = \{g_i, g_{i+1}\}$ and $R^{\mathbf{B}}(\{c, d\}, \mathbf{p}_j) = \{g_j, g_{j+1}\}$. Then we have $R^{\mathbf{B}}(c, \mathbf{p}_j) = g_j$ or $R^{\mathbf{B}}(d, \mathbf{p}_j) = g_j$. Since $g_j \wedge g_i = g_j = g_j \wedge g_{i+1}$, we have

$$g_j = R^{\mathbf{B}}(c, \mathbf{p}_j) \wedge R^{\mathbf{B}}(d, \mathbf{p}_i) \quad \text{or} \quad g_j = R^{\mathbf{B}}(d, \mathbf{p}_j) \wedge R^{\mathbf{B}}(c, \mathbf{p}_i).$$

By the compatibility of the operation $\wedge$, either conclusion leads to

$$g_j = R^{\mathbf{B}}(c \wedge d, \mathbf{p}_i \wedge \mathbf{p}_j).$$

But we also have $R^{\mathbf{B}}(c, \mathbf{p}_j) = g_{j+1}$ or $R^{\mathbf{B}}(d, \mathbf{p}_j) = g_{j+1}$. The same argument yields that
$$g_{j+1} = R^{\mathbf{B}}(c \wedge d, \mathbf{p}_i \wedge \mathbf{p}_j),$$
which contradicts $g_j > g_{j+1}$.

We have shown that, deleting trivial links in the descending chain $g_1, \ldots, g_n$ if necessary, the polynomials associated to the links may be chosen to be unary polynomials constructed from members of $P$, and that the polynomials associated with different links may be assumed to be constructed from different members of $P$. This shows that the length of the descending chain $g_1, \ldots, g_n$, after trivial links are deleted, is at most $|P|$. We get a similar bound on the ascending chain $h_1, \ldots, h_n$. But, there are only finitely many special principal congruence formulas which describe Mal'cev chains of length $\leq 2|P|$ and involve only terms in $P$. Let $\Phi$ be the disjunction of all such special principal congruence formulas. $\Phi$ is a principal congruence formula with the property that for a randomly chosen instance of $(a, b) \in \text{Cg}^{\mathbf{B}}(c, d)$ we have $\Phi^{\mathbf{B}}(a, b, c, d)$. Hence $\Phi$ is a principal congruence formula for $\mathcal{V}$. $\square$

This theorem proves that a finitely generated variety which has a compatible semilattice term operation has DPC. It is not true that every locally finite variety with a compatible semilattice term operation has DPC as the next example shows.

**Example 4.4** Our algebra will be $\mathbf{A} = (A; \wedge, \cdot)$ where $A = \{0, 1, 2, 3, \ldots\}$ and
$$x \wedge y = \begin{cases} x & \text{if } x = y \\ 0 & \text{otherwise} \end{cases}, \quad x \cdot y = \begin{cases} x & \text{if } y = x + 1 \\ 0 & \text{otherwise} \end{cases}.$$
Here $+$ is the usual addition on $\{0, 1, 2, 3, \ldots\}$. The operation $\wedge$ is a semilattice operation. It is easy to check that $(x \cdot y) \wedge (u \cdot v)$ is nonzero if and only if $0 \neq x = u = y - 1 = v - 1$ in which case the value of the expression is $x$. The same statement is true of $(x \wedge u) \cdot (y \wedge v)$. This implies that
$$(x \cdot y) \wedge (u \cdot v) = (x \wedge u) \cdot (y \wedge v),$$
or that $\wedge$ is a compatible semilattice operation of $\mathbf{A}$.

For any $x, y \in A$ we have $x \wedge y, x \cdot y \in \{x, y, 0\}$, so for any $S \subseteq A$ the subuniverse generated by $S$ is contained in $S \cup \{0\}$. Hence $\mathbf{A}$ is a uniformly locally finite algebra which implies that $\mathcal{V} := \mathcal{V}(\mathbf{A})$ is locally finite. However, $\mathcal{V}$ does not have DPC. The easiest way to establish this is for us to refer to Theorem 5.1 to see that $\mathcal{V}$ is residually small. Next, we refer to the result of [1] which proves that any residually small variety *which has DPC* is residually $< N$ for some finite $N$. Finally, $\mathcal{V}$ is not residually $< N$ for any finite $N$ since $\mathbf{A}$ is an infinite subdirectly irreducible algebra in $\mathcal{V}$.

# 5 Residual Smallness and Finite Axiomatizability

McKenzie has characterized which locally finite varieties with type–set contained in $\{\mathbf{3},\mathbf{4},\mathbf{5}\}$ are residually small. One can apply McKenzie's criterion to show that any locally finite variety which has a compatible semilattice term is residually small. If $\mathbf{A}$ is a finite algebra with a compatible semilattice term operation, then McKenzie's result implies that $\mathcal{V}(\mathbf{A})$ is residually $< (2^\omega)^+$. However, by the previously mentioned result from [1], this forces $\mathcal{V}(\mathbf{A})$ to be residually $< N$ for some finite $N$. Without using McKenzie's results (which are not yet published) we prove that if $\mathbf{A}$ has a compatible semilattice term operation and $|A| = n$, then $\mathcal{V}(\mathbf{A})$ is residually less than $N = 2^{n^{n^{n^{n+1}}}}$.

**THEOREM 5.1** *Let $\mathcal{V}$ be a variety which has a compatible semilattice term operation. $\mathcal{V}$ is residually small. If $\mathcal{V} = \mathcal{V}(\mathbf{A})$ where $|A| = n$, then $\mathcal{V}$ is residually less than $2^{n^{n^{n^{n+1}}}}$.*

*Proof.* If $\mathcal{V}$ is not finitely generated, then let $P$ be a minimal set of $\mathcal{V}$–inequivalent representatives for the terms of $\mathcal{V}$. If $\mathcal{V}$ is generated by $\mathbf{A}$ where $|A| = n$, then let $P$ be a minimal set of representatives for the $(n^n + 1)$–ary terms. In either case, a typical member of $P$ will be written $t(x, \mathbf{y})$.

Let $\mathbf{B}$ be a subdirectly irreducible algebra in $\mathcal{V}$. Whether or not $\mathcal{V}$ is finitely generated, for each $p(x) \in \mathrm{Pol}_1(\mathbf{B})$ there is a $t(x,\mathbf{y}) \in P$ and a tuple $\mathbf{b}$ of elements of $B$ such that $p(x) = t^{\mathbf{B}}(x, \mathbf{b})$. (When $\mathcal{V}$ is finitely generated this claim follows from Corollary 4.2.) If $\mu$ is the monolith of $\mathbf{B}$, choose elements $0 < 1$ in $B$ such that $\mu = \mathrm{Cg}(0,1)$. For each $u \neq v$ in $B$ we get that $(0,1) \in \mathrm{Cg}(u,v)$, so there is a Mal'cev chain connecting 1 to 0 by a chain of polynomial images of $\{u,v\}$. We may assume that the form of this Mal'cev chain is of the kind considered when we constructed our DPC formula in the last section, so in particular we may assume that this Mal'cev chain is a chain descending from 1 followed by a chain ascending to 0. Hence whenever $u \neq v$ in $B$ there is a $t(x, \mathbf{y}) \in P$ and a tuple $\mathbf{b}$ such that

$$t^{\mathbf{B}}(u, \mathbf{b}) = 1 > t^{\mathbf{B}}(v, \mathbf{b})$$

or the same with $u$ and $v$ switched. (This situation occurs at the beginning of the Mal'cev chain.)

For each $u \in B$, let $S_u$ denote the subset of $P$ consisting of those $t(x,\mathbf{y}) \in P$ for which there is a tuple $\mathbf{b}$ of elements of $B$ such that $t^{\mathbf{B}}(u, \mathbf{b}) = 1$. Assume for the moment that $u \neq v$ and that $S_u = S_v$. Since $u \neq v$, we have from the previous paragraph that there is some $t \in P$ and some tuple $\mathbf{b}$ such that

$$t^{\mathbf{B}}(u, \mathbf{b}) = 1 > t^{\mathbf{B}}(v, \mathbf{b})$$

or the same with $u$ and $v$ switched. Without loss of generality we may assume that the displayed situation is our case. The displayed line implies that $t \in$

$S_u = S_v$, so there is some tuple $\mathbf{b}'$ such that $t^{\mathbf{B}}(v, \mathbf{b}') = 1$. Hence

$$\begin{aligned}
t^{\mathbf{B}}(v, \mathbf{b}) &< 1 \\
&= 1 \wedge 1 \\
&= t^{\mathbf{B}}(u, \mathbf{b}) \wedge t^{\mathbf{B}}(v, \mathbf{b}') \\
&= t^{\mathbf{B}}(u \wedge v, \mathbf{b} \wedge \mathbf{b}') \\
&= t^{\mathbf{B}}(u \wedge v, \mathbf{b}' \wedge \mathbf{b}) \\
&= t^{\mathbf{B}}(u, \mathbf{b}') \wedge t^{\mathbf{B}}(v, \mathbf{b}) \\
&\leq t^{\mathbf{B}}(v, \mathbf{b})
\end{aligned}$$

But $t^{\mathbf{B}}(v, \mathbf{b}) < t^{\mathbf{B}}(v, \mathbf{b})$ is impossible. We conclude that if $u \neq v$, then $S_u \neq S_v$. Consequently, the function from $B$ to the power set of $P$ defined by $u \mapsto S_u$ is one–to–one. This proves that $|B| \leq 2^{|P|}$, and therefore that $\mathcal{V}$ is residually less than $(2^{|P|})^+$. When $\mathcal{V}$ is generated by an $n$–element algebra, then $|P| = |F_{\mathcal{V}}(n^n + 1)|$. The estimate $|F_{\mathcal{V}}(m)| \leq n^{n^m}$, which holds for any algebra in a variety generated by an $n$–element algebra, yields that

$$|B| \leq 2^{n^{n^{(n^n+1)}}} \leq 2^{n^{n^{n^{n+1}}}}.$$

The latter inequality is strict if $n > 1$ and the former is strict if $n = 1$. This proves the theorem. $\square$

**COROLLARY 5.2** *A finite algebra of finite type which has a compatible semilattice term operation is finitely based.*

*Proof.* Any finite algebra of finite type which generates a variety with DPC which is residually $< N$ for some integer $N$ is finitely based, as is proved in [5]. $\square$

# 6 Idempotent Varieties

The self–rectangulating conditon is a kind of term condition similar to the usual one which defines the class of abelian algebras. There seem to be several points of correspondence between the property of being abelian and the property of being self–rectangulating. In this correspondence one term condition is the analogue of the other.

We find a point of correspondence when we restrict the types that can occur in a finite algebra. A finite abelian algebra must have type–set contained in $\{\mathbf{1}, \mathbf{2}\}$. When the algebra is of type–set $\{\mathbf{2}\}$ only, then it satisfies a property stronger than the term condition: it is quasiaffine. This means that it is a subalgebra of a reduct of an algebra which has a compatible Mal'cev term operation. Correspondingly, a finite self–rectangulating algebra of type–set $\{\mathbf{5}\}$ is a subalgebra of a reduct of an algebra which has a compatible semilattice term operation.

Another point of correspondence is found when we consider locally finite varieties with the properties described in the previous paragraph. An abelian locally finite variety $\mathcal{V}$ of type–set $\{\mathbf{2}\}$ is congruence modular, and therefore affine by Herrmann's Theorem. But a variety is affine if and only if it has a compatible Mal'cev term operation. Thus a locally finite variety of type–set $\{\mathbf{2}\}$ is abelian if and only if it has a compatible Mal'cev term operation. Correspondingly, we have proved in this paper that a locally finite variety of type–set $\{\mathbf{5}\}$ is self–rectangulating if and only if it has a compatible semilattice term operation.

The correspondence up to this point suggests that the analogue of an affine algebra is an algebra with a compatible semilattice operation. Now the structure of a variety $\mathcal{V}$ of affine algebras is understood by associating to the variety a ring $\mathbf{R}$ and then comparing $\mathcal{V}$ to the variety of left $\mathbf{R}$–modules. The ring of $\mathcal{V}$ is constructed from the compatible Mal'cev term operation of $\mathcal{V}$ along with some idempotent binary terms. In particular, the ring associated to $\mathcal{V}$ is the same as the ring associated to the idempotent reduct of $\mathcal{V}$. It is satisfying to learn that the correspondence between abelian varieties and self–rectangulating varieties continues to hold up. We will see in this section that if $\mathcal{V}$ has a compatible semilattice term operation, then we can construct a *semiring* from the semilattice term operation and some idempotent binary term operations. When $\mathcal{V}$ is idempotent, the semiring seems to determine almost as much about $\mathcal{V}$ as the ring of an idempotent affine variety determines about that variety.

The results in this section come from [3]. The paper [3] is part of an investigation of varieties of idempotent algebras whose term operations commute with each other. That paper is an analysis of the case where the variety has a semilattice term operation. Hence all varieties considered in [3] have a compatible semilattice term operation, but moreover they have the property that all other term operations are compatible as well. If one checks the proofs given, though, one finds that almost all arguments remain intact (sometimes with inessential modifications) for idempotent varieties with a compatible semilattice operation. We state the main facts about the structure of such varieties. (N. B.: in [3] the compatible semilattice operation is a join operation denoted by $+$. In this paper we have used a compatible meet operation denoted by $\wedge$. Therefore, the order we have been considering is the reverse of the order in [3]. Since we mean to emphasize the analogy with affine algebras, and therefore with varieties of modules, we switch now to a join semilattice operation denoted by $+$.)

**THEOREM 6.1** (The Structure of Subdirectly Irreducible Algebras) *Let $\mathbf{A}$ be an idempotent subdirectly irreducible algebra which has a compatible (join) semilattice term operation denoted by $+$. The following are true.*

(1) *$\mathbf{A}$ has a least element $0$.*

(2) *$\mathbf{A}$ has an element $u$ which is the least element in $A - \{0\}$.*

(3) The monolith of **A** is the equivalence relation on $A$ having $\{0, u\}$ as its only nontrivial block.

(4) The set
$$\begin{aligned} U &= \{s(x,0) \mid s \text{ is a binary term}\} \\ &= \{p(x) \in \text{Pol}_1(\mathbf{A}) \mid p(0) = 0\} \\ &= \{q(x) \in \text{Pol}_1(\mathbf{A}) \mid 0 \in q(A)\} \end{aligned}$$
is a set of decreasing endomorphisms of $(A; +)$.

(5) **A** is polynomially equivalent to $(A; +, U)$.

*Proof.* See the proof of Theorem 3.1 of [3]. □

Now we construct the semiring. We could build it from a set of inequivalent binary terms or from the elements of $\mathbf{F}_\mathcal{V}(x, y)$ represented by those binary terms. We follow [3] and choose the latter approach.

Let $\mathcal{V}$ be an idempotent variety with a compatible semilattice term operation denoted by $+$. Let $R$ be the subuniverse of $\mathbf{F}_\mathcal{V}(x, y)$ consisting of all $t \in F_\mathcal{V}(x, y)$ such that $t + y = t$. We write $e_t$ for the endomorphism of $\mathbf{F}_\mathcal{V}(x, y)$ determined by $x \mapsto t, y \mapsto y$. For $s, t \in F_\mathcal{V}(x, y)$ we write $s \circ t$ to denote $e_t(s)$. We write $0$ to denote the element $y \in F_\mathcal{V}(x, y)$ and we write $1$ to denote $x + y \in F_\mathcal{V}(x, y)$.

**Definition 6.2** $\mathbf{R}(\mathcal{V})$ is the algebra of type $\langle 2, 2, 0, 0 \rangle$ given by $(R; \circ, +, 1, 0)$. $\mathbf{R}(\mathcal{V})$ is called the **semiring of the variety** $\mathcal{V}$.

$\mathbf{R}(\mathcal{V})$ satisfies the associative laws with respect to $\circ$ and $+$ and $\circ$ distributes over $+$ on either side. We have $0 \circ x = x \circ 0 = 0$ while $1 \circ x = x \circ 1 = x$. The operation $+$ is commutative and $0 + x = x$. (Cf. the proof of Theorem 4.11 of [3].) These laws define what is called a "semiring" in [3]. However, $\mathbf{R}(\mathcal{V})$ satisfies an additional equation not valid in every semiring: $1 + x = 1$.

Next we define the **coefficient representation of the term operations of** $\mathcal{V}$. Let $t = t(x_1, \ldots, x_n)$ be an arbitrary $n$–ary term operation of $\mathcal{V}$. Since $\mathcal{V}$ is idempotent and $+$ is a compatible semilattice term operation of $\mathcal{V}$, we have that
$$\mathcal{V} \models t(x_1, \ldots, x_n) + y = \sum_{i=1}^{n} (t(y, y, \ldots, y, x_i, y, \ldots, y) + y).$$

(On the right hand side $x_i$ is in the $i$–th position.) This shows that $t$ is uniquely determined by the elements $\widehat{t}_1, \ldots, \widehat{t}_n$ of $\mathbf{R}(\mathcal{V})$ where $\widehat{t}_i$ denotes the element of $R \subseteq F_\mathcal{V}(x, y)$ equal to $t(y, y, \ldots, y, x, y, \ldots, y) + y$ with $x$ in the $i$–th position only. The coefficient representation of $t$ is
$$\widehat{t}_1 \bullet x_1 + \cdots + \widehat{t}_n \bullet x_n,$$
which may be thought of as a term in the language of $\mathbf{R}(\mathcal{V})$–semimodules. (See [3] for the definition of an $\mathbf{R}(\mathcal{V})$–semimodule.)

**THEOREM 6.3** (The Connection with Semimodules) *The assignment*

$$t(x_1, \ldots, x_n) \mapsto \widehat{t}_1 \bullet x_1 + \cdots + \widehat{t}_n \bullet x_n$$

*is an injective homomorphism from the clone of $\mathcal{V}$ to the clone of $\mathbf{R}(\mathcal{V})$–semimodules.*

*Proof.* This follows from Lemma 3.5 of [3]. □

The homomorphism of the previous theorem is not surjective since the clone of $\mathbf{R}(\mathcal{V})$–semimodules is not idempotent. However, even if we replace the target clone with the clone of all idempotent $\mathbf{R}(\mathcal{V})$–semimodule operations then the image of the homomorphism still may not be surjective. (It is in some cases, for example when 1 is join–irreducible in $\mathbf{R}(\mathcal{V})$. See the first claim in Theorem 4.24 of [3]. We note, however, that for a variety $\mathcal{V}$ generated by a single subdirectly irreducible algebra, 1 is not necessarily join–irreducible in $\mathbf{R}(\mathcal{V})$ if $\mathbf{R}(\mathcal{V})$ is not commutative.) Nevertheless, the connection with semimodules described in Theorem 6.3 is useful for proving results about $\mathcal{V}$. It can be used to prove the following theorem.

**THEOREM 6.4** (The Lattice of Equational Theories) *The lattice of equational theories extending the equational theory of $\mathcal{V}$ is isomorphic (in a natural way) to $\mathrm{Con}\,(\mathbf{R}(\mathcal{V}))$.*

*Proof.* See the proof of Theorem 4.20 of [3]. □

Any module variety has an injective cogenerator. A weaker statement is true for idempotent varieties with a compatible semilattice term operation. One does not have an injective cogenerator, usually, but there is a "canonical" cogenerator constructible from $\mathbf{R}(\mathcal{V})$ whose structure is very comprehensible. We begin the construction of this cogenerating algebra now.

**Definition 6.5** Let $\mathbf{R}$ be a semiring satisfying $1 + x = 1$. An **annihilator ideal** of $\mathbf{R}$ is a subset $I \subseteq R$ which is an order ideal and is closed under $+$.

It can be shown that the annihilator ideals of $\mathbf{R}$ are precisely the subsets of $\mathbf{R}$ of the form $0/\theta$ where $\theta$ is a congruence on $\mathbf{R}$.

Let $\mathcal{I}$ denote the set of annihilator ideals of $\mathbf{R}(\mathcal{V})$. For $a, b \in \mathcal{I}$ let $a \oplus b$ denote $a \cap b$. For each $r \in R(\mathcal{V})$ and each $a \in \mathcal{I}$ let

$$(a)r^{-1} = \{s \in R(\mathcal{V}) \mid sr \in a\}.$$

If $f$ is a $n$–ary basic operation symbol of $\mathcal{V}$ and the coefficient representation of $f$ is

$$f(x_1, \ldots, x_n) = \widehat{f}_1 \bullet x_1 + \cdots + \widehat{f}_n \bullet x_n,$$

then let $[f]$ denote the $n$–ary operation on $\mathcal{I}$ defined by

$$[f](a_1, \ldots, a_n) = (a_1)\widehat{f}_1^{-1} \oplus \cdots \oplus (a_n)\widehat{f}_n^{-1}.$$

We let $\mathcal{I}(\mathcal{V})$ denote the algebra with universe $\mathcal{I}$ and with basic operations of the form $[f]$, one for each basic operation of $\mathcal{V}$.

**THEOREM 6.6** (The Canonical Cogenerator) *With notation as above, the following hold.*

(1) $\mathcal{I}(\mathcal{V}) \in \mathcal{V}$.

(2) *If* $\mathbf{A} \in \mathcal{V}$ *is a subdirectly irreducible algebra, then* $\mathbf{A}$ *is embeddable into* $\mathcal{I}(\mathcal{V})$.

(3) $\mathcal{V} = \mathbf{SP}(\mathcal{I}(\mathcal{V}))$.

*Proof.* This can be proved with the arguments in Subsection 4.3 of [3]. Some care is needed to put things on the left or right, as needed, since $\mathbf{R}(\mathcal{V})$ may not be commutative. □

Our last result cited from [3] concerns the congruence extension property. One may view this result as a strengthening of Theorem 4.3 for idempotent varieties, since any locally finite variety with the congruence extension property has DPC (see [1]).

**THEOREM 6.7** (The Congruence Extension Property) *An idempotent variety with a compatible semilattice term operation has the congruence extension property.*

*Proof.* The proof in Section 5 of [3] works here. □

Assume that $\mathcal{V}$ is an affine variety and that $t(x_1, \ldots, x_n)$ is a term operation of $\mathcal{V}$. Let $t'(x) = t(x, x, \ldots, x)$. Then the term operation

$$s(x_1, x_2, \ldots, x_n) = t(x_1, x_2, \ldots, x_n) - t'(x_1) + x_1$$

is an idempotent term operation and $t$ can be reconstructed from $s$ and $t'$ since

$$t(x_1, x_2, \ldots, x_n) = s(x_1, x_2, \ldots, x_n) - x_1 + t'(x_1).$$

Thus the clone of $\mathcal{V}$ is generated by its idempotent subclone and the unary component of its clone. This is part of the reason why the ring associated to $\mathcal{V}$ determines most of the information about $\mathcal{V}$. The corresponding result for varieties with a compatible semilattice term operation is false. The idempotent subclone may be only a very small part of the clone. A consequence of this is that the analogy between self–rectangulating varieties of type–set $\{\mathbf{5}\}$ and affine varieties breaks down (or seems to) when one tries to pass from idempotent varieties to varieties which are not idempotent. For example, in the affine case, a variety with a finite associated ring can only have finitely many subvarieties. For *idempotent* varieties with a compatible semilattice operation the corresponding result is true: if the associated semiring is finite then there are only finitely many subvarieties. But for nonidempotent varieties with a compatible semilattice operation one can have a finite associated semiring and still have infinitely many subvarieties. The algebra of Example 4.4 generates a variety of this kind.

We remark that idempotence is necessary to obtain all the results of this section. For varieties with a compatible semilattice term which are not idempotent none of the results in this section are true. For each result stated in this section, either Example 4.4 or the modification of Example 2.10 obtained by adjoining a compatible semilattice operation is a counterexample to the result when idempotence is omitted as a hypothesis.

Department of Mathematical Sciences, University of Arkansas, Fayetteville, AR 72701, USA.

Bolyai Institute, Aradi vértanúk tere 1, H–6720 Szeged, Hungary.